\documentclass[11pt,a4paper]{article}
\usepackage{amssymb}
\usepackage{amsmath}
\usepackage[numbers,sort]{natbib}
\usepackage{hyperref}
\usepackage{epsfig}
\usepackage{graphics}
\usepackage{caption}
\usepackage{subcaption}
\makeatletter
\newcommand{\labeleditem}[2]{%
\def\@itemlabel{\textbf{#1}}
\item
\def\@currentlabel{#1}\label{#2}}
\makeatother

\newtheorem{theorem}{Theorem}[section]

\newtheorem{lemma}[theorem]{Lemma}

\newcommand{\Ac}{\mathcal{A}}
\newcommand{\Aca}{\mathcal{A}_{\alpha}}
\newcommand{\Acah}{\mathcal{A}_{\alpha,h}}
\newcommand{\Lc}{\mathcal{L}}
\newcommand{\Bc}{\mathcal{B}}
\newcommand{\Bca}{\mathcal{B}_{\alpha}}

\newcommand{\Mc}{\mathcal{M}}
\newcommand{\Vc}{\mathcal{V}}

\newcommand{\MINRES}{\textsc{Minres}}

\newcommand{\Ltwo}{{L^2(\Omega)}}
\newcommand{\LtwoB}{{L^2(\partial \Omega)}}
\newcommand{\Hone}{{H^1(\Omega)}}
\newcommand{\Htwo}{{H^2(\Omega)}}
\newcommand{\bHtwo}{{\bar{H}^2(\Omega)}}

\newcommand{\aLtwo}{L_{\alpha}^2(\Omega)}
\newcommand{\aiLtwo}{L_{\alpha^{-1}}^2(\Omega)}

\newcommand{\aHtwo}{H_{\alpha}^2(\Omega)}

\makeatletter
\renewcommand{\inf}{\mathop{\@inf\vphantom{\@sup}}}
\renewcommand{\sup}{\mathop{\@sup\vphantom{\@inf}}}
\newcommand{\@inf}{\operatorname*{inf}}
\newcommand{\@sup}{\operatorname*{sup}}
\makeatother

\newcommand{\hess}{\nabla^2}     % Hessian operator
\newcommand{\grad}{\nabla}       % Gradient
\newcommand{\laplace}{\Delta}    % Laplacian

\renewcommand{\ker}{\operatorname{Ker}} % Kernel
\newcommand{\cond}{\kappa}              % Condition number

\newcommand{\nvec}{\mathbf{n}} % Normal vector
\newcommand{\ze}{{0}}           % Matrix zero entry
\newcommand{\subs}{\,\cdot\,}  % Substitution/placeholder symbol
\newcommand{\dx}[1][x]{\, d#1} % integration dx

\makeatletter
\newcommand{\norm}{\@ifstar\@lrnorms\@norm}
\newcommand{\@lrnorm}[1]{%
  \left\lVert
   #1
  \right\rVert
}
\newcommand{\@norm}[2][]{% 
  \mathopen{#1\lVert}
  #2
  \mathclose{#1\rVert}
}
\newcommand{\abs}{\@ifstar\@lrabs\@abs}
\newcommand{\@lrabs}[1]{% 
  \left\lvert
   #1
  \right\rvert
}
\newcommand{\@abs}[2][]{%
  \mathopen{#1\lvert}
  #2
  \mathclose{#1\rvert}
}
\renewcommand{\brack}{\@ifstar\@lrbrack\@brack}
\newcommand{\@lrbrack}[1]{%
  \left\langle
   #1
  \right\rangle
}
\newcommand{\@brack}[2][]{%
  \mathopen{#1\langle}
  #2
  \mathclose{#1\rangle}
}
\newcommand{\pp}{\@ifstar\@lrpp\@pp}
\newcommand{\@lrpp}[1]{%
  \left(
   #1
  \right)
}
\newcommand{\@pp}[2][]{%
  \mathopen{#1(}
  #2
  \mathclose{#1)}
}
\newcommand{\set}{\@ifstar\@lrset\@set}
\newcommand{\@lrset}[1]{%
  \left\{
   #1
  \right\}
}
\newcommand{\@set}[2][]{%
  \mathopen{#1\{}
  #2
  \mathclose{#1\}}
}
\makeatother

\title{Robust preconditioners for PDE-constrained optimization with limited observations} 
\author{Kent-Andre Mardal\thanks{{\em Kent-Andre Mardal}, Center for Biomedical Computing, Simula Research Laboratory; Department of Mathematics, University of Oslo, Norway. Email: {\tt kent-and@simula.no}}, Bj{\o}rn Fredrik Nielsen\thanks{{\em Bj{\o}rn Fredrik Nielsen}, Department of Mathematical Sciences and Technology,
Norwegian University of Life Sciences, Norway; Simula Research Laboratory;  
Center for Cardiological Innovation, Oslo University Hospital. 
Email: {\tt bjorn.f.nielsen@nmbu.no}} and Magne Nordaas\thanks{{\em Magne Nordaas}, Center for Biomedical Computing, Simula Research Laboratory. Email: {\tt magneano@simula.no}}}

\begin{document}
\maketitle 

\begin{abstract}
  Regularization robust preconditioners for PDE-constrained
  optimization problems have been successfully developed. These
  methods, however, typically assume that observation data is
  available throughout the entire domain of the state equation. For
  many inverse problems, this is an unrealistic assumption.  In this
  paper we propose and analyze preconditioners for PDE-constrained
  optimization problems with limited observation data,
  e.g. observations are only available at the boundary of the solution
  domain. Our methods are robust with respect to both the
  regularization parameter and the mesh size. That is, the condition
  number of the preconditioned optimality system is uniformly bounded, independently
  of the size of these two parameters. We first consider a
  prototypical elliptic control problem and thereafter more general
  PDE-constrained optimization problems. Our theoretical findings are
  illuminated by several numerical results.
\end{abstract}



\newpage
\section{Introduction}
Consider the model problem: 
\begin{equation}
\label{A1}
\min_{f, \, u} \,
\set*{ 
  \frac{1}{2} \norm{u-d}^2_{\LtwoB)} 
  + \frac{\alpha}{2} \norm{f}_{\Ltwo}^2
},
\end{equation}
on a Lipschitz domain $\Omega\subset \mathbb{R}^n$, subject to
\begin{alignat}{2}
\label{A2}
- \laplace u + u + f &= 0 
\quad &&\mbox{in } \Omega, \\
\label{A3} 
\frac{\partial u}{\partial \mathbf{n}} &= 0 
\quad &&\mbox{on } \partial \Omega. 
\end{alignat}
This minimization task is similar to the standard example considered in PDE-constrained 
optimization. But instead of assuming that observation data is available everywhere in $\Omega$, 
we consider the case where observations are only given at the boundary $\partial \Omega$ of 
$\Omega$, that is $d \in \LtwoB$, see the first term in \eqref{A1}. 
For problems of the form \eqref{A1}-\eqref{A3}, in which 
\begin{equation}
  \label{A4}
  \frac{1}{2}\norm{u-d}^2_{\LtwoB}  + \frac{\alpha}{2} \norm{f}_{\Ltwo}^2
\end{equation}
is replaced by
\begin{equation}
  \frac{1}{2}\norm{u-d}^2_{\Ltwo}  + \frac{\alpha}{2} \norm{f}_{\Ltwo}^2
\end{equation}
very efficient preconditioners have been developed for the associated
KKT system. In fact, by employing proper $\alpha$-dependent scalings
of the involved Hilbert spaces \cite{s-z}, or by using a Schur
complement approach \cite{Pea12}, methods that are robust with respect
to the size of the regularization parameter $\alpha$ have been
developed. 
More specifically, the condition number of the preconditioned optimality system is
small and bounded independently of $0 < \alpha \ll 1$
and the mesh size $h$. This ensures good performance for suitable
Krylov subspace methods, e.g. the minimum residual method (\MINRES{}),
independently of both parameters.
%More specifically, the convergence rate of suitable
%preconditioned Krylov subspace methods (e.g. \MINRES{}) applied to the 
%discretized optimality system can be estimated independently of $0 <
%\alpha \ll 1$ and the mesh size $h$.
These techniques have been
extended to handle time dependent problems \cite{Pea12_II} and
PDE-constrained optimization with Stokes equations \cite{z-11}, but
the rigorous analysis of $\alpha$-independent bounds always requires
that observations are available throughout all of $\Omega$.

For cases with limited observations, for example with cost-functionals
of the form \eqref{A4}, efficient preconditioners are also available for a
rather large class of PDE-constrained optimization problems, see
\cite{Nie13,Nie10}. But these techniques do not yield convergence rates, for
the preconditioned KKT-system, that are completely robust with respect
to the size of the regularization parameter $\alpha$. Instead, the number
of preconditioned \MINRES{} iterations grows logarithmically\footnote{In
  \cite{Nie13,Nie10} it is proved that the number of needed preconditioned
  \MINRES{} iterations cannot grow faster than \[a + b \left[ \log_{10}
    \left( \alpha^{-1} \right) \right]^2.\] Furthermore, in \cite{Nie13}
  it is explained why iterations counts of the kind \eqref{A6} often
  will occur in practice.} with respect to the size of $\alpha^{-1}$, as $\alpha \rightarrow
0$:
\begin{equation}
\label{A6}
a + b \log_{10} \left( \alpha^{-1} \right). 
\end{equation}
According to the numerical experiments presented in \cite{Nie13}, the size of $b$ may become 
significant. More specifically, $b \in [5,50]$ for problems with simple elliptic state equations 
posed on rectangles. Thus, for small values of $\alpha$, \MINRES{} may require rather many iterations 
to converge - even though the growth in iteration numbers is only logarithmically. 

In practice, observations are rarely available throughout the entire
domain of the state equation.  On the contrary, the purpose of solving
an inverse problem is typically to use data recorded at the surface of
an object to compute internal properties of that object: Impedance
tomography, the inverse problem of electrocardiography (ECG),
computerized tomography (CT), etc. This fact, combined with the
discussion above, motivate the need for further improving numerical
methods for solving KKT systems arising in connection with
PDE-constrained optimization.

This paper is organized as follows.  In the next section we derive the
KKT system associated with the model problem \eqref{A1}-\eqref{A3}.  Our
$\alpha$ robust preconditioner is presented in Section
\ref{Numerical_experiments}, along with a number of numerical
experiments. Sections \ref{Analysis}-\ref{Preconditioning} contain our
analysis, and the method is generalized in Sections
\ref{Generalization}-\ref{Alternative}. Section \ref{Conclusions}
provides a discussion of our findings, including their limitations.

\section{KKT system} 
\label{KKT_system}
Consider the PDE \eqref{A2} with the boundary condition \eqref{A3}.  A
solution $u$ to this elliptic PDE, with source term $f\in\Ltwo$, is
known to have improved regularity, i.e.  $u\in H^{1+s}(\Omega)$, for
some $s\in [0, 1]$, with $s$ depending on the domain $\Omega$. In the remainder
of this paper we assume that $u$ has full regularity, i.e. $u\in
\Htwo$.  This is known to hold if $\Omega$ is convex or if
$\partial\Omega$ is $C^2$, see e.g. \cite{grisvard1985, BHac92}.

When solutions to \eqref{A2} exhibit this improved regularity, we can
write the problem on the non-standard variational form: Find $u \in \bHtwo$ such
that
\begin{equation}
\label{B1} 
(- \laplace u + u,w)_{\Ltwo} + (f,w)_{\Ltwo} = 0\quad \forall w \in \Ltwo,  
\end{equation}
where 
\begin{equation*}
\bHtwo =  \bigg\{
\phi \in \Htwo \, \bigg| \, 
\frac{\partial \phi}{\partial \mathbf{n}} = 0  
\mbox{ on } \partial \Omega
\bigg\}, 
\end{equation*}
equipped with the inner product
\begin{equation}
  \label{eq:H2_inner}
  \begin{split}
      (u,v)_{\Htwo} 
      &= \int_\Omega \hess u : \hess v + \grad u \cdot \grad v + uv \dx \\
      &= \int_\Omega \laplace u \laplace v + \grad u \cdot \grad v + u v \dx.
  \end{split}
\end{equation}
Here $\hess u$ denotes the Hessian of $u$, and the second identity is due to the boundary condition $\frac{\partial u}{\partial \nvec} = 0$
imposed on the space $\bHtwo$.

We will see below that, in order to design a regularization robust
preconditioner for \eqref{A1}-\eqref{A3}, it is convenient to express the
state equation in the form \eqref{B1}, instead of employing integration
by parts/Green's formula to write it on the standard self-adjoint form. 
%If $f \in \Ltwo$ and $\partial \Omega$ is sufficiently regular/smooth, then 
%standard theory for elliptic PDEs assures that \eqref{B1} has a unique solution 
%$u \in \bHtwo$.   

\subsection{Optimality system} 
We may express \eqref{A1}-\eqref{A3} in the form: 
\begin{equation}
\label{B3}
\min_{f \in \Ltwo, \, u \in \bHtwo}\,
\set*{ 
  \frac{1}{2}\norm{u-d}_{\LtwoB}^2 + \frac{\alpha}{2} \norm{f}_{\Ltwo}^2
}
\end{equation}
subject to 
\begin{equation}
\label{B4} 
(- \laplace u + u,w)_{\Ltwo} + (f,w)_{\Ltwo}  =0\quad \forall w \in \Ltwo.   
\end{equation}
The associated Lagrangian reads 
\begin{equation*}
\begin{split}
  \mathcal{L}(f,u,w) &= 
  \frac{1}{2} \norm{u-d}_{\LtwoB}^2 
  + \frac{\alpha}{2} \norm{f}_{\Ltwo}^2
  +(f - \Delta u + u,w)_{\Ltwo}, 
\end{split}
\end{equation*}
with $f \in \Ltwo$, $u \in \bHtwo$ and $w \in \Ltwo$. 
From the first order optimality conditions 
\begin{equation*} 
 \frac{\partial \mathcal{L}}{\partial f} = 0,\quad
 \frac{\partial \mathcal{L}}{\partial u} = 0,\quad 
 \frac{\partial \mathcal{L}}{\partial w} = 0, 
\end{equation*} 
we obtain the optimality system: Determine 
$(f,u,w) \in \Ltwo \times \bHtwo \times \Ltwo$ such that  
\begin{alignat}{2}
  \label{B5.1}
  \alpha (f,\psi)_{\Ltwo}+(\psi,w)_{\Ltwo} &= 0 &&\quad \forall \psi \in \Ltwo, \\
  \label{B5.2}
  (u-d,\phi)_{\LtwoB} + (-\laplace \phi + \phi, w)_{\Ltwo} &= 0 &&\quad \forall \phi \in \bHtwo, \\
  \label{B5.3}
  (f,\xi)_{\Ltwo} + (- \laplace u + u,\xi)_{\Ltwo}  &= 0 &&\quad \forall \xi \in \Ltwo.
\end{alignat} 

\section{Numerical experiments}
\label{Numerical_experiments} 
Prior to analyzing our model problem, we will consider some numerical experiments. 
Discretization of \eqref{B5.1}-\eqref{B5.3} yields an algebraic system of the form 
\begin{equation}
\label{N1}
\underbrace{\left[
    \begin{array}{ccc}
      \alpha M & \ze        & M   \\
      \ze      & M_{\partial} & A^T \\
      M        & A          & 0 
    \end{array}
  \right]}_{\Aca}
\left[
  \begin{array}{c}
    f \\ u \\ w
  \end{array}
\right]
=
\left[
  \begin{array}{c}
    0           \\ 
    \tilde M_{\partial}d \\ 
    0
  \end{array}
\right], 
\end{equation}
where 
\begin{itemize}
\item $M$ is a mass matrix, 
\item $M_{\partial}$ is a mass matrix associated with the boundary $\partial \Omega$ of $\Omega$. 
\item $A$ is a matrix that arise upon discretization of
  the operator $(1-\Delta)$. Since we write the state equation on a non self-adjoint
  form, $A$ will not be the usual sum of the stiffness and mass
  matrices. Instead, equation \eqref{B1} is discretized with subspaces
  of $\bHtwo$ and $\Ltwo$.   
\end{itemize}

In the current numerical experiments, we employ the Bogner-Fox-Schmit
(BFS) rectangle for discretizing the state variable $u \in \bHtwo$.
That is, the finite element field consists of bicubic polynomials
that are continuous, have continuous first order derivatives and mixed
second order derivatives at each vertex of the mesh.  BFS elements are
$C^1$ on rectangles and therefore $H^2$-conforming. The control $f$
and Lagrange multiplier $w$ are discretized with discontinuous bicubic
elements.

We propose to precondition \eqref{N1} with the block-diagonal matrix
\begin{equation}
\label{N2}
\Bca=
\left[
\begin{array}{ccc}
\alpha M & \ze               & \ze \\
\ze      & \alpha R+M_\partial & \ze \\
\ze      &  \ze              & \frac{1}{\alpha} M 
\end{array}
\right]^{-1}, 
\end{equation} 
%% For the FEM spaces $V_h$ and $Q_h$, it turns out that 
%% $A^TM^{-1}A=B$, 
where $R$ results from a discretization of the bilinear form $b(\cdot,\cdot)$ on $\bHtwo$: 
\begin{equation}
  b( u, v)
%%   &=& \int_\Omega D^2 u : D^2 v + 2 Du\cdot Dv + uv \, dx \\
\label{N2.1}
   = (u,v)_{H^2(\Omega)} + \int_\Omega \grad u \cdot \grad v \, dx.  
\end{equation}
In the experiments presented below, we used this bilinear form to
construct a multigrid approximation of $\left( \alpha R+M_\partial
\right)^{-1}$.

\subsection*{Remark} 
The bilinear form \eqref{N2.1} is equivalent to the inner product on
$\bHtwo$. The additional term stems from our choice of implementing a
multigrid algorithm for the bilinear form associated with the operator
$(\laplace - 1)^2 = \laplace^2 -2 \laplace +1$.  Indeed, the bilinear
form $\alpha b(\subs, \subs) + (\subs,\subs)_{\LtwoB}$ can be seen to
coincide with the variational form associated with the fourth order problem
\begin{alignat*}{2}
\label{A2}
\alpha(\laplace-1)^2 u  &= f 
\quad &&\mbox{in } \Omega, \\
\frac{\partial u}{\partial \mathbf{n}} &= 0 
\quad &&\mbox{on } \partial \Omega, \\
\alpha\frac{\partial \laplace u}{\partial \mathbf{n}} &= u 
\quad &&\mbox{on } \partial \Omega. 
\end{alignat*}
%% solver for
%% \begin{eqnarray*}
%% (\Delta - I)^2 \phi &=& g \quad \mbox{in } \Omega, \\
%% \frac{\partial \phi}{\partial \mathbf{n}} &=& 0 
%% \quad \mbox{on } \partial \Omega,
%% \end{eqnarray*}   
%% in order to realize an isomorphism between $\bHtwo$ and $\bHtwo'$.

\subsection*{}
\begin{figure}
      \centering
      \begin{subfigure}[t]{.46\textwidth}
        \includegraphics[width=1.0\textwidth]{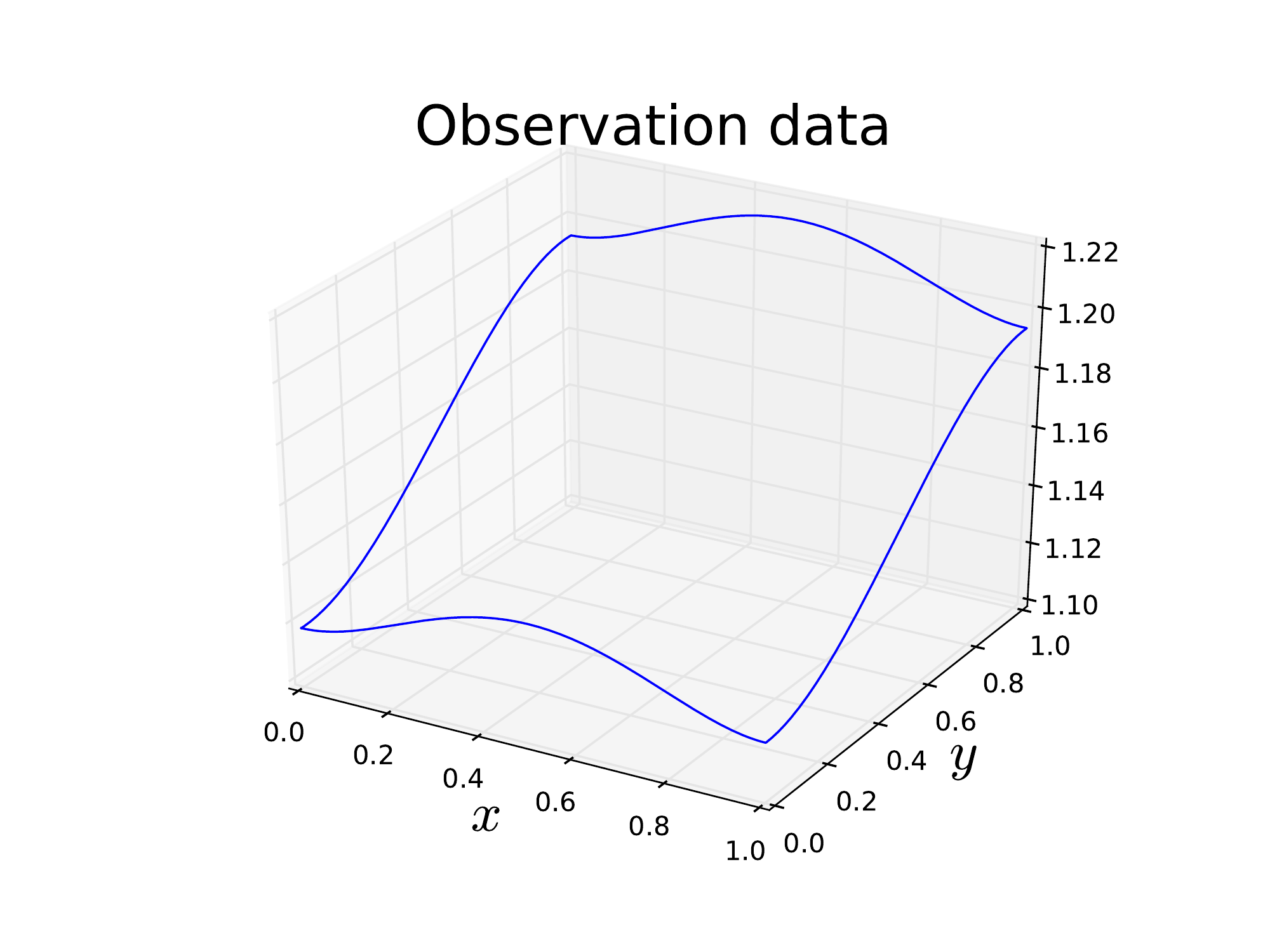}
        \captionsetup{font=scriptsize}
        \caption{Observation data $d$. The forward model was solved for
          the control shown in (\subref{fig2:d}), but only the
          boundary values can be observed.}
        \label{fig2:b}
      \end{subfigure}
            \hfill
      \centering
      \begin{subfigure}[t]{.46\textwidth}
        \includegraphics[width=1.0\textwidth]{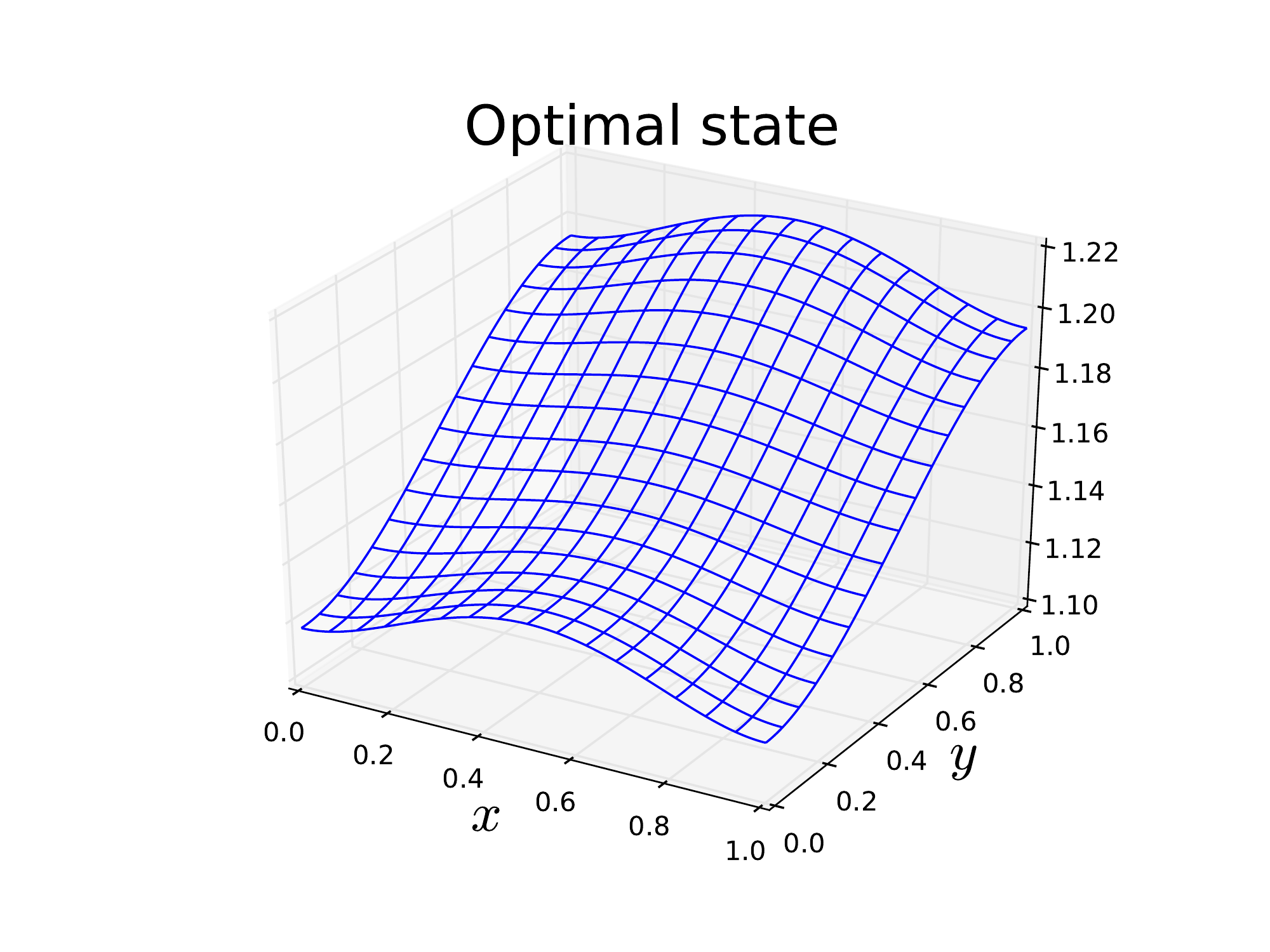}
        \captionsetup{font=scriptsize}
        \caption{Computed optimal state $u$ based on the observation data
          shown in (\subref{fig2:b})}
        \label{fig2:a}
      \end{subfigure}
      \centering
      \begin{subfigure}[t]{.46\textwidth}
        \includegraphics[width=1.0\textwidth]{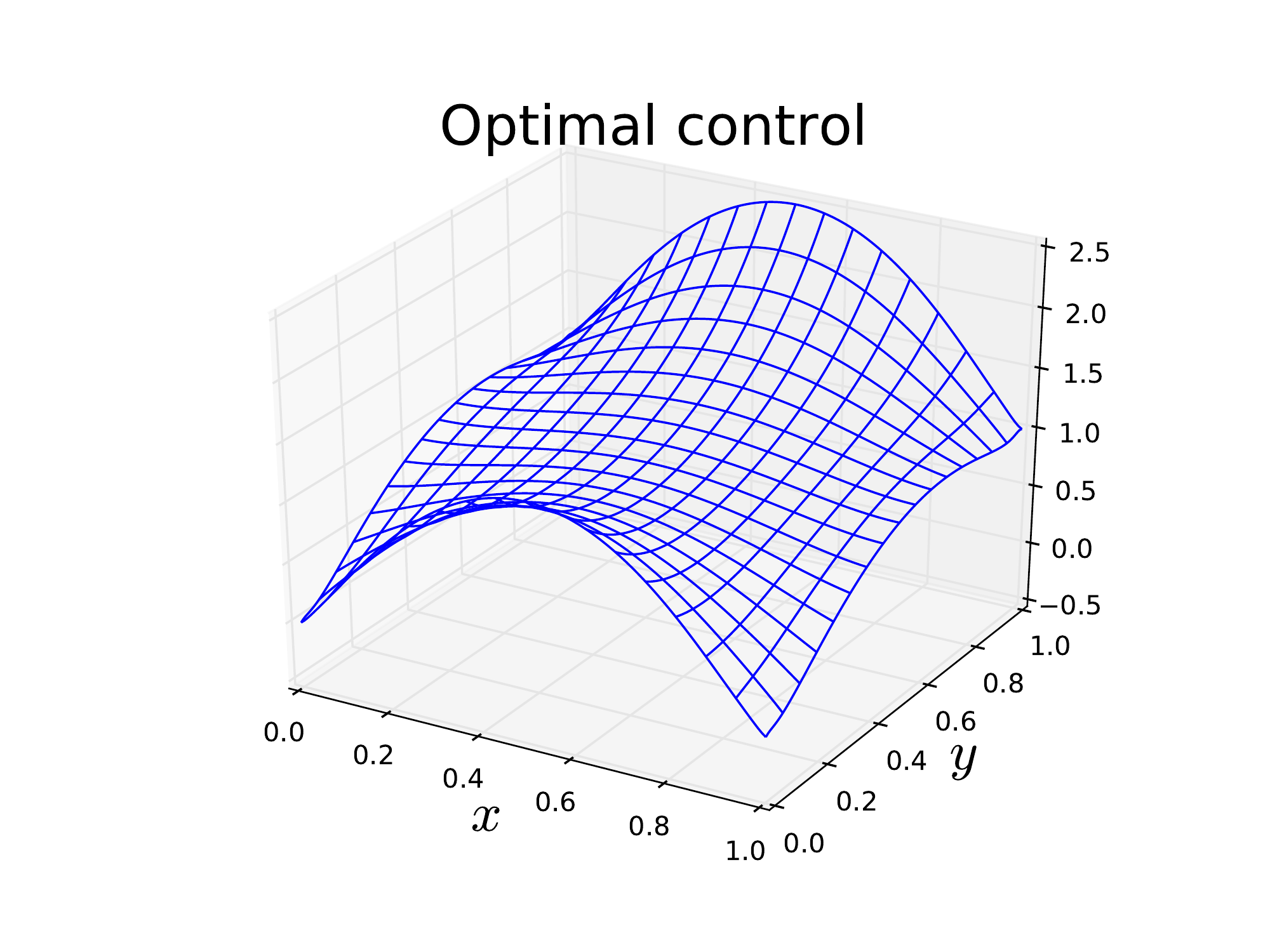}
        \captionsetup{font=scriptsize}
        \caption{Computed optimal control $f$ based on the observation
          data in (\subref{fig2:b})}
        \label{fig2:c}
      \end{subfigure}
      \hfill
      \centering
      \begin{subfigure}[t]{.46\textwidth}
        \includegraphics[width=1.0\textwidth]{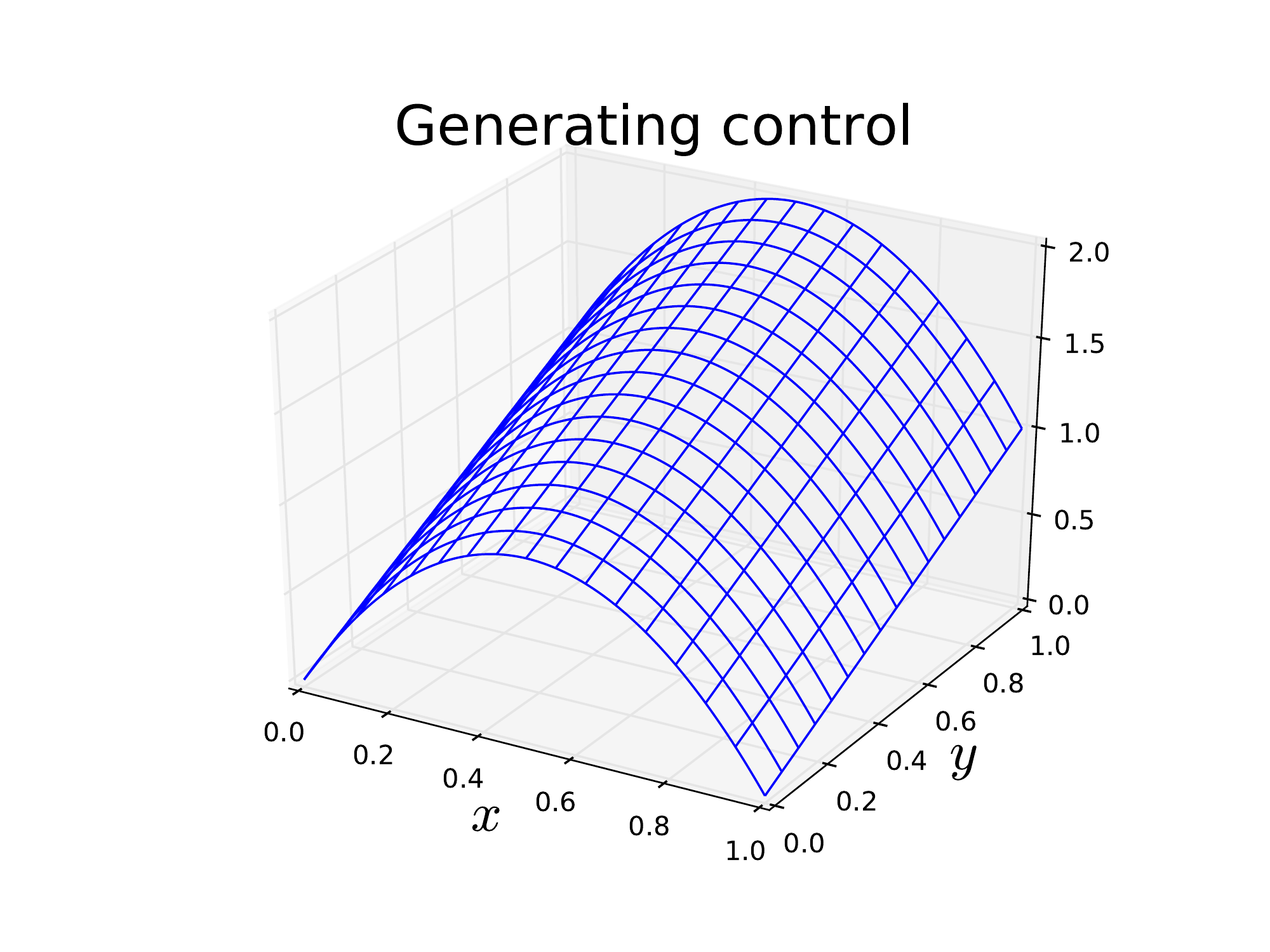}
        \captionsetup{font=scriptsize}
        \caption{The ``true'' control function used to generate the observation data in
          (\subref{fig2:b}).}
        \label{fig2:d}
      \end{subfigure}
      \caption{An example of a solution of \eqref{N1}.  The
        observation data $d$ was generated with the forward model, using 
        the ``true'' control $4x(1-x) + y$ shown in panel 
        (\subref{fig2:d}). Solutions to the unregularized problem are
        non-unique, and the generating control cannot be (exactly) recovered. The
        figures were generated with mesh parameter $h = 1/128$
        and regularization parameter $\alpha = 10^{-6}$.}
      \label{fig:2}
\end{figure}

To limit the technical complexity of the implementation, we considered
the problem \eqref{A1}-\eqref{A3} on the unit square in two dimensions.
The experiments were implemented in Python and SciPy. The meshes were
uniform rectangular, with the coarsest level for the multigrid solver
consisting of $8\times 8$ rectangles. Figure \ref{fig:2} shows an example of a 
solution of the optimality system \eqref{N1}. 
%Further details about the discretization is presented in Appendix \ref{Discretization}.

\subsection{Eigenvalues} 
Let us first consider the exact preconditioner $\Bca$ defined in \eqref{N2}. 
If $\Bca$ is a good preconditioner for the
discrete optimality system \eqref{N1}, then the spectral condition
number of $\Bca \Aca$ should be small and bounded, independently
of the size of both the regularization parameter $\alpha$ and the discretization
parameter $h$. 

\begin{figure}
  \centering
      \centering
      \begin{subfigure}[b]{.46\textwidth}
        \includegraphics[width=1.0\textwidth]{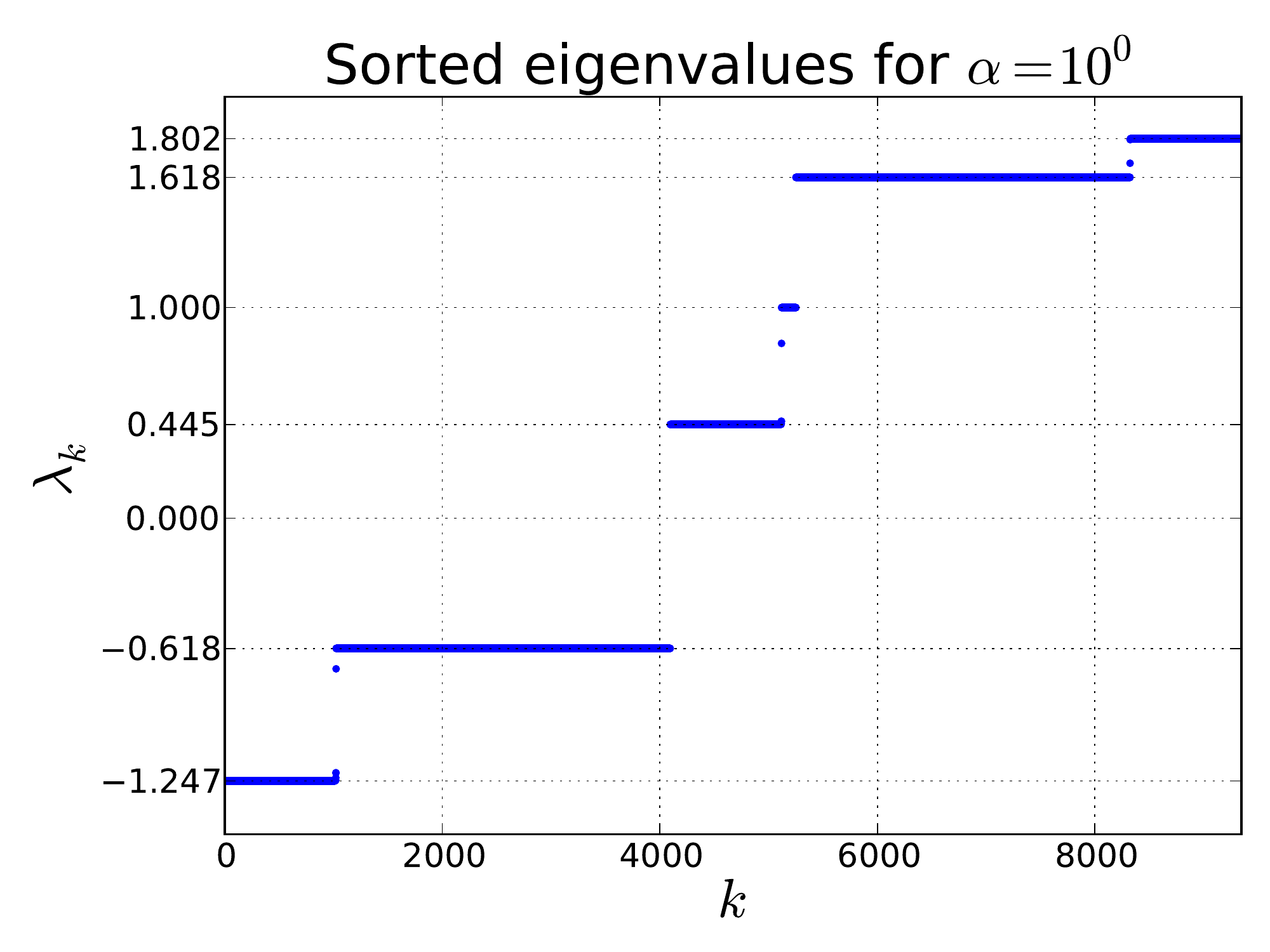}
        \captionsetup{font=scriptsize}
        \label{fig:1a}
        \caption{$\alpha = 1$}
      \end{subfigure}
      \hfill
      \centering
      \begin{subfigure}[b]{.46\textwidth}
        \includegraphics[width=1.0\textwidth]{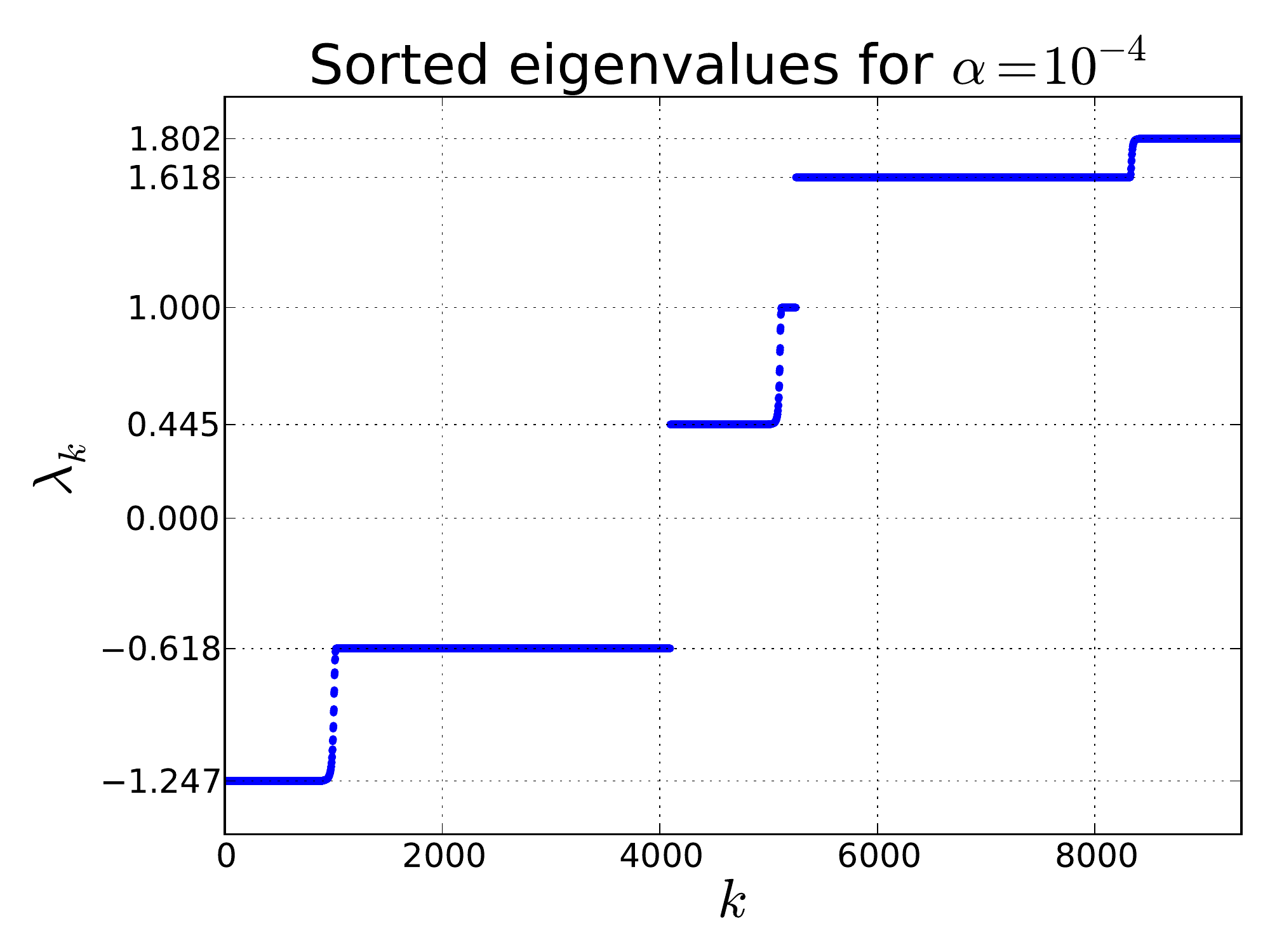}
        \captionsetup{font=scriptsize}
        \label{fig:1b}
        \caption{$\alpha = 10^{-4}$}
      \end{subfigure}
      \centering
      \begin{subfigure}[b]{.46\textwidth}
        \includegraphics[width=1.0\textwidth]{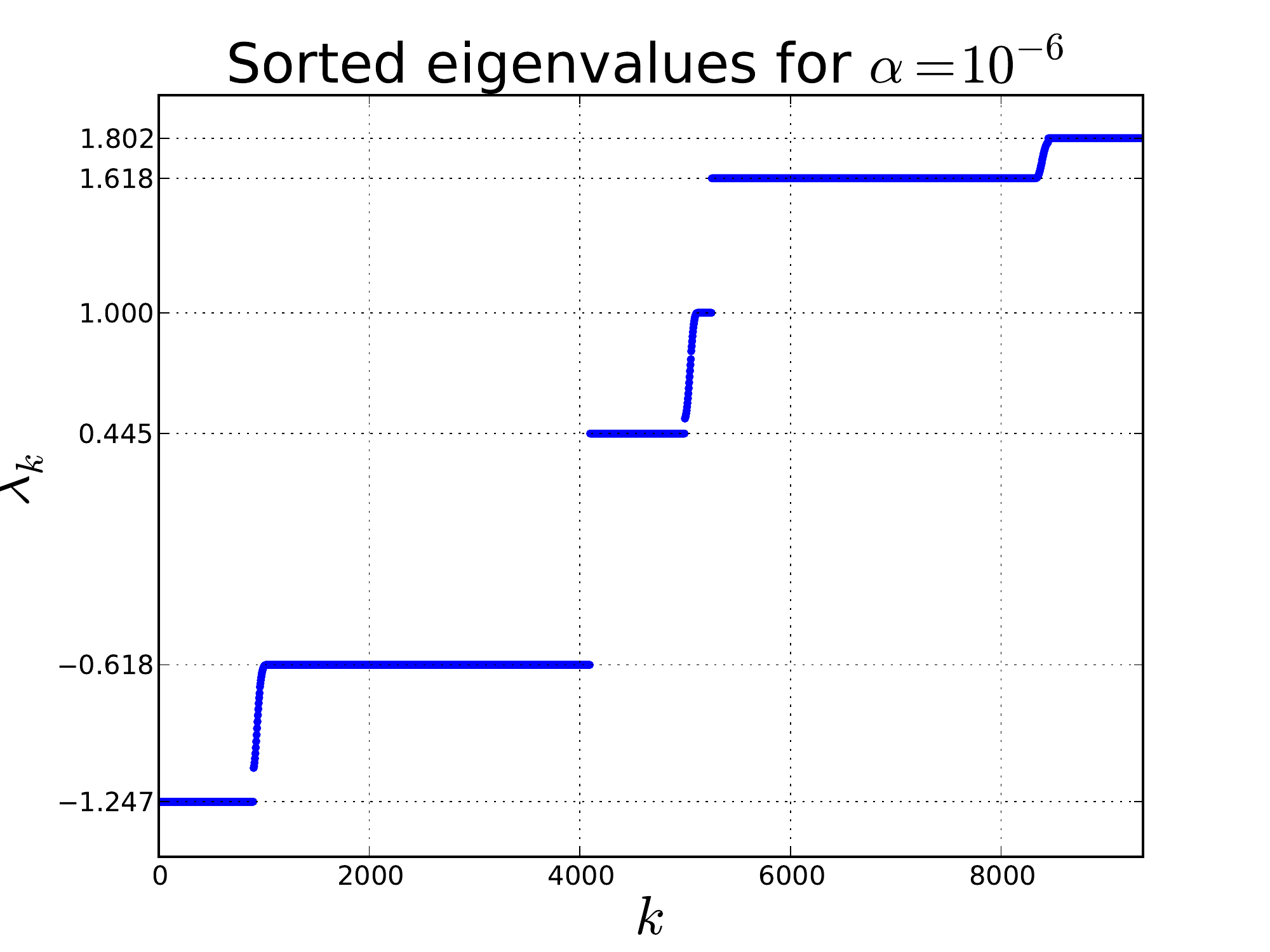}
        \captionsetup{font=scriptsize}
        \label{fig:1c}
        \caption{$\alpha = 10^{-6}$}
      \end{subfigure}
      \hfill
      \centering
      \begin{subfigure}[b]{.46\textwidth}
        \includegraphics[width=1.0\textwidth]{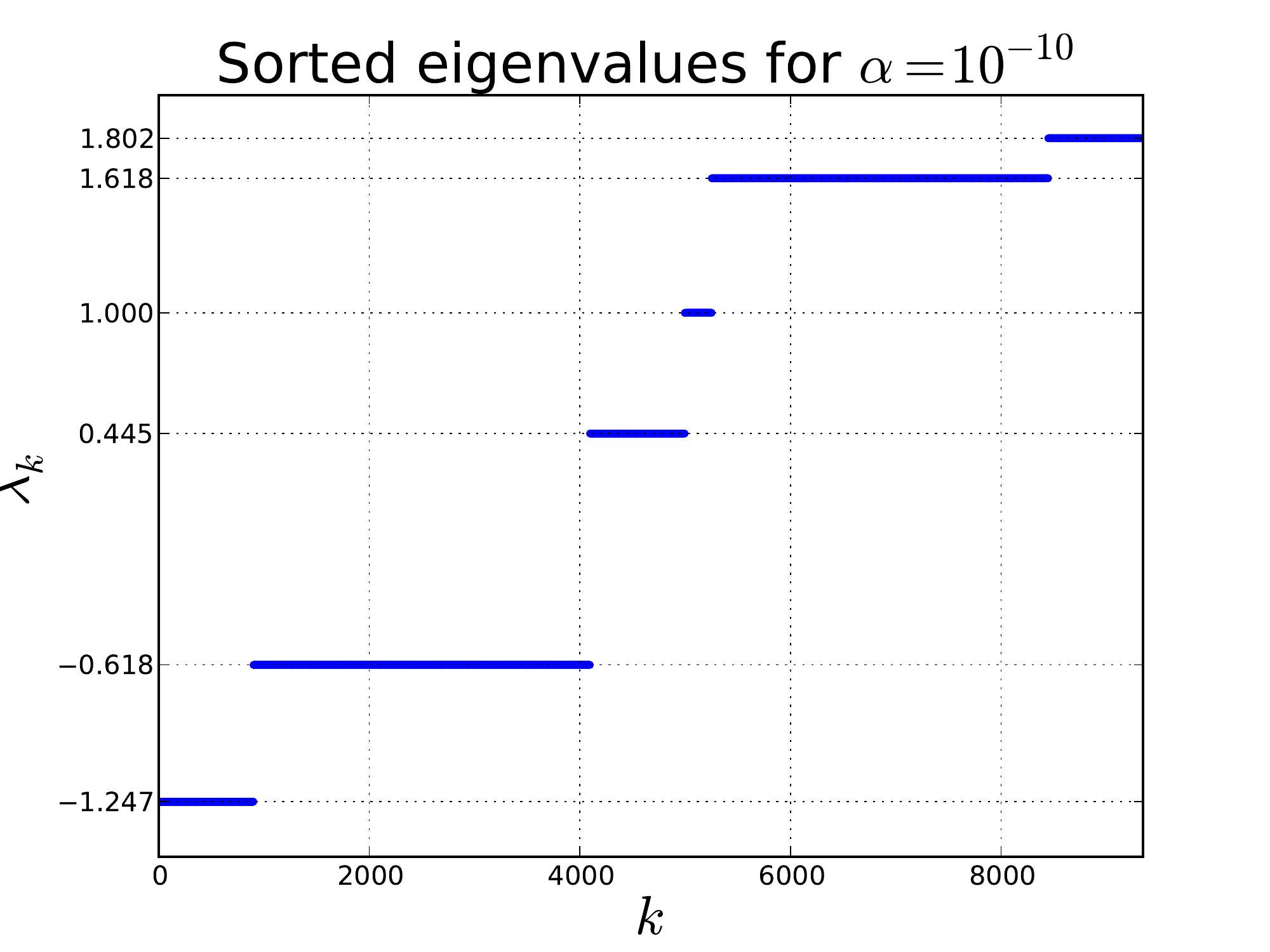}
        \captionsetup{font=scriptsize}
        \label{fig:1d}
        \caption{$\alpha = 10^{-10}$}
      \end{subfigure}
  \caption{Spectrum of $\Bca \Aca$ for different regularization
    parameters $\alpha$. The discretization parameter
    was $h=2^{-4}$ for all figures.}
  \label{fig:1}
\end{figure}

The eigenvalues of this preconditioned system were
computed by solving the generalized eigenvalue problem
\begin{equation*}
  \Aca x = \lambda \Bca^{-1} x.
\end{equation*}
We found that the absolute value of the eigenvalues $\lambda$ were bounded, with 
\begin{equation*}
  0.445 \leq \vert\lambda\vert \leq 1.809,
\end{equation*}
uniformly in $\alpha\in\{1,10^{-1},\ldots, 10^{-10}\}$ and $h \in
\{2^{-2},\ldots,2^{-5}\}$. This yields a uniform condition number
$k(\Bca \Aca) \approx 4.05$. The spectra of the preconditioned
systems are pictured in Figure \ref{fig:1} for some
choices of $\alpha$. The spectra are clearly divided into three
bounded intervals, and the eigenvalues are more clustered for $\alpha
\approx 1$ and for very small $\alpha$.

\subsection{Efficient preconditioning}
In practice, the action of $\Bca$ is replaced with a less
computationally expensive operation $\widehat{\Bca}$. Note that $\Bca$
has a block structure, and that computationally efficient approximations can be
constructed for the individual blocks.  Specifically, $\widehat{\Bca}$
is constructed by employing
\begin{itemize}
\item 1 multigrid V-cycle  for the (2,2) block of $\Bca$, containing a symmetric $4\times 4$ block Gauss-Seidel 
smoother where the blocks contain the matrix entries corresponding to  all degrees of freedom associate with a vertex in the mesh
(see \cite{Zhang1994} for a theoretical analysis of the method).   
\item 2 symmetric Gauss-Seidel iterations for the (1,1) and (3,3) blocks.
\end{itemize}
%% The multigrid cycle for bicubic BFS discretization was analyzed in
%% \cite{Zhang1994}. We remark that the BFS element results in nested
%% discrete spaces, in the sense that $V_h\supset V_H$ if
%% $\mathcal{T}_h$ is a refinement of $\mathcal{T}_H$.  This property,
%% while ideal for multigrid, does not hold for the majority of the
%% developed $C^1$-conforming elements.  Multigrid methods for non-nested
%% but conforming spaces was analyzed in \cite{BrambleZhang1995}.

We estimated condition numbers of the individual blocks of $\Bca^{-1}$ preconditioned
with their respective approximations. The results are reported in Tables
\ref{tab:gs} and \ref{tab:multigrid}.  A slight deterioration in the
performance of the multigrid cycle can be seen for very small values of 
$\alpha > 0$. 
\begin{table}%
  \label{tab:block_conditon_numbers}
  \centering
  \begin{minipage}{.45\textwidth}%
    \begin{tabular}{c|ccc}
      Iterations  & 1     & 2     & 3    \\
      \hline
       ($h=2^{-8}$)& 1.931 & 1.303 & 1.126
    \end{tabular}
    \caption{Condition numbers of $M$ preconditioned with symmetric
      Gauss-Seidel iterations.}
    \label{tab:gs}
    \qquad
  \end{minipage}%
  \qquad
  \begin{minipage}{.45\textwidth}%
    \begin{tabular}{l|ccc}
      $\alpha$\textbackslash$h$ 
               & $2^{-4}$ & $2^{-6}$ & $2^{-8}$ \\
      \hline
      $1$      &  1.130  & 1.136   & 1.140   \\
      $10^{-4}$ &  1.129  & 1.135   & 1.139   \\
      $10^{-8}$ &  1.237  & 1.150   & 1.149   \\
      $10^{-12}$&  1.252  & 1.259   & 1.253      
    \end{tabular}
    \caption{Estimated condition numbers of $\alpha R + M_\partial$ 
      preconditioned with one V-cycle multigrid iteration.}
    \label{tab:multigrid}
  \end{minipage}%
\end{table}

\subsection{Iteration numbers}
To verify that also $\widehat{\Bca}$ is an effective
preconditioner for $\Aca$, we applied the \MINRES{} scheme to the system
\begin{equation*}
  \widehat\Bca \Aca x = \widehat\Bca b.  
\end{equation*}
For the results presented in Table
\ref{tab:MINRES_iterations_and_estimates}, the \MINRES{} iteration
process was stopped as soon as
\begin{equation}
\label{N3}
\frac{(r_k,\widehat{\Bca} r_k)}{(r_0,\widehat{\Bca} r_0)}= 
\frac{(\Aca x_k - b, \widehat{\Bca} \{\Aca x_k - b\})}{(\Aca x_0 - b, \widehat{\Bca} \{\Aca x_0 - b\})} 
\leq \varepsilon, 
\end{equation}
which is the standard termination criterion for the preconditioned
\MINRES{} scheme, provided that the preconditioner is SPD. A random
initial guess $x_0$ was used, and the tolerance was set to $\varepsilon =
10^{-12}$. 

\begin{table}
  \centering
  \begin{tabular}{l|ccccccccccc}
    $\alpha$\textbackslash$h$ & $ 2^{-4}$&$2^{-5}$& $2^{-6}$& $2^{-7}$ \\ 
    \hline
     1        & 53(4.33) & 53(4.36) & 53(4.36) & 53(4.36) \\
     $10^{-1}$ & 57(4.31) & 57(4.34) & 57(4.35) & 57(4.35) \\
     $10^{-2}$ & 75(4.31) & 72(4.34) & 70(4.35) & 68(4.35) \\
     $10^{-3}$ & 79(4.31) & 79(4.34) & 77(4.35) & 73(4.35) \\
     $10^{-4}$ & 81(4.30) & 81(4.33) & 79(4.35) & 77(4.35) \\
     $10^{-5}$ & 82(4.33) & 81(4.33) & 79(4.35) & 79(4.35) \\
     $10^{-6}$ & 81(4.35) & 79(4.36) & 79(4.35) & 81(4.35) \\
     $10^{-7}$ & 70(4.35) & 81(4.37) & 81(4.36) & 79(4.35) \\
     $10^{-8}$ & 62(4.36) & 70(4.36) & 79(4.36) & 81(4.36) \\
     $10^{-9}$ & 62(4.36) & 64(4.37) & 68(4.37) & 78(4.36) \\
     $10^{-10}$& 62(4.36) & 63(4.36) & 64(4.37) & 67(4.37) 
  \end{tabular}
  \caption{Number of preconditioned \MINRES{} iterations needed to solve the optimality system to a relative error tolerance $\varepsilon = 10^{-12}$. 
    Estimated condition numbers in parentheses, computed from conjugate gradient iterations on the normal equations for the preconditioned optimality system.}% $\widehat{\Bca} \Aca$.}
  \label{tab:MINRES_iterations_and_estimates}
\end{table}

\section{Analysis of the KKT system}
\label{Analysis}
Recall that our optimality system reads: 
\begin{alignat*}{2}
\alpha (f,\psi)_{\Ltwo}+(\psi,w)_{\Ltwo} &= 0 &&\quad \forall \psi \in \Ltwo, \\
(u-d,\phi)_{\LtwoB} + (- \Delta \phi + \phi, w)_{\Ltwo} &= 0 &&\quad \forall \phi \in \bHtwo, \\
(f,\xi)_{\Ltwo} + (- \Delta u + u,\xi)_{\Ltwo}  &= 0 &&\quad \forall \xi \in \Ltwo, 
\end{alignat*} 
with unknowns $f \in \Ltwo$, $u \in \bHtwo$ and $w \in \Ltwo$. 
We may write this KKT system in the form: \\ Determine 
$(f,u,w) \in \Ltwo \times \bHtwo \times \Ltwo$ such that 
\begin{equation}
\label{B5}
\underbrace{\left[
\begin{array}{ccc}
\alpha M & \ze        & M' \\
     \ze & M_{\partial} & A' \\
       M &  A         & 0 
\end{array}
\right]}_{\Aca}
\left[
\begin{array}{c}
f \\ u \\ w
\end{array}
\right]
=
\left[
\begin{array}{c}
0 \\ 
\tilde{M}_{\partial}d \\ 
0
\end{array}
\right], 
\end{equation}
where 
\begin{align}
\label{B6}
 M&: \Ltwo \rightarrow \Ltwo', &  f &\mapsto (f,\subs)_{\Ltwo}, \\ 
\label{B7}
 M_{\partial}&: \bHtwo \rightarrow \bHtwo', &  u &\mapsto (u,\subs)_{\LtwoB}, \\ 
\label{B8}
 \tilde{M}_{\partial} &: \LtwoB \rightarrow \bHtwo', &  d &\mapsto (d,\subs)_{\LtwoB}, \\ 
\label{B9}
 A &: \bHtwo \rightarrow \Ltwo', & u &\mapsto (- \Delta u + u,\subs)_{\Ltwo},  
\end{align}
and the notation "$'$" is used to denote dual operators and dual spaces. 
In the rest of this paper, the symbols $M$, $ M_{\partial}$ and $A$ will represent the mappings 
defined in \eqref{B6}, \eqref{B7} and \eqref{B9}, respectively, and not (the associated) matrices, as 
was the case in Section \ref{Numerical_experiments}. (We believe that this mild ambiguity improves the readability of the present text). 

By using standard techniques for saddle point problems, one can show
that the system \eqref{B5} satisfies the Brezzi conditions~\cite{brezzi},
provided that $\alpha >0$.  Therefore, for every $\alpha > 0$, this
set of equations has a unique solution.  Nevertheless, if the standard
norms of $\Ltwo$ and $\Htwo$ are employed in the analysis, then the
constants in the Brezzi conditions will depend on $\alpha$. More
specifically, the constant in the coercivity condition will be of
order $O(\alpha)$, and thus becomes very small for $0 < \alpha \ll
1$. This property is consistent with the ill posed nature of
\eqref{A1}-\eqref{A3} for $\alpha=0$, and makes it difficult to design
$\alpha$ robust preconditioners for the algebraic system associated
with \eqref{B5}. 

Similar to the approach used in \cite{s-z, Mar11, Nie10}, we will now introduce
weighted Hilbert spaces. The weights are constructed such that the
constants appearing in the Brezzi conditions are independent of
$\alpha$. Thereafter, in Section \ref{Preconditioning}, we will show
how these scaled Hilbert spaces can be combined with simple maps to
design $\alpha$ robust preconditioners for our model problem.

\subsection{Weighted norms}
Consider the $\alpha$-weighted norms:  
\begin{align}
\label{B0.1}
 \norm{ f }_{\aLtwo}^2 &= \alpha \norm{f}_{\Ltwo}^2, \\
\label{B0.2}
 \norm{ u }_{\aHtwo}^2 &= \alpha \norm{ u }_{\Htwo}^2 + \norm{u}_{\LtwoB}^2, \\
\label{B0.3}
 \norm{ w }_{\aiLtwo}^2 &= \frac{1}{\alpha} \norm{ w }_{\Ltwo}^2, 
\end{align} 
applied to the control $f$, the state $u$ and the dual/Lagrange-multiplier $w$, respectively. Note that these norms become
``meaningless'' for $\alpha = 0$, but are well defined for positive $\alpha$.

\subsection{Brezzi conditions}
We will now analyze the properties of 
\begin{equation*}
  \Aca: \aLtwo \times \aHtwo \times \aiLtwo \rightarrow \aLtwo' \times \aHtwo' \times \aiLtwo',    
\end{equation*}
defined in \eqref{B5}. 
%% where 
%% \[
%% \Aca = 
%% \left[
%% \begin{array}{ccc}
%% \alpha M & 0 & M' \\
%% 0 & M_{\partial} & A' \\
%% M &  A & 0 
%% \end{array}
%% \right],
%% \]
%% and $M$, $M_{\partial}$ and $A$ are defined in \eqref{B6}, \eqref{B7} and \eqref{B9}, respectively. 
More specifically, we will show that the Brezzi conditions are satisfied 
with constants that do not depend on the size of the regularization parameter $\alpha > 0$. 
Note that we use the scaled Hilbert norms \eqref{B0.1}-\eqref{B0.3}. 

%% If these weighted norms are employed, then we can establish suitable
%% bounds with constants that are independent of the size of the
%% regularization parameter $\alpha$. 
\begin{lemma} 
\label{pre_inf-sup}
For all  $\alpha > 0$, the following ``inf-sup'' condition holds:
\begin{equation*}
\inf_{w \in \aiLtwo} \sup_{(f,u) \in \aLtwo \times \aHtwo} 
\frac{(f,w)_{\Ltwo} + (- \laplace u + u,w)_{\Ltwo}}{\norm{ (f,u) }_{\aLtwo \times \aHtwo} \norm{ w }_{\aiLtwo}} \geq 1.
\end{equation*}
\end{lemma}
\subsubsection*{Proof} 
Note that $\aLtwo$ and $\aiLtwo$ contain the same functions, provided that $\alpha > 0$. 
Let $w \in \aiLtwo$ be arbitrary. 
By choosing $f=w$ and $u=0$ we find that 
\begin{equation*} 
  \begin{aligned}
    \sup_{(f,u) \in \aLtwo \times \aHtwo} 
    \frac{(f,w)_{\Ltwo} + (- \Delta u + u,w)_{\Ltwo}}{\norm{ (f,u) }_{\aLtwo \times \aHtwo} \norm{ w }_{\aiLtwo}} &\geq
    \frac{(w,w)_{\Ltwo}}{\norm{ (w,0) }_{\aLtwo \times \aHtwo}  \norm{ w }_{\aiLtwo}} \\ 
    &= 
    \frac{\norm{ w }_{\Ltwo}^2}{\sqrt{\alpha}\norm{ w }_{\Ltwo} (\sqrt{\alpha})^{-1} \norm{ w }_{\Ltwo}} \\
    &= 
    1.
  \end{aligned}
\end{equation*}
Since $w\in \aiLtwo$ was arbitrary, this completes the proof.
\\
\rule{2mm}{2mm} \\
Expressed in terms of the operators that constitute 
$\Aca$, Lemma \ref{pre_inf-sup} takes the form 
\begin{equation*}
\inf_{w \in \aiLtwo} \sup_{(f,u) \in \aLtwo \times \aHtwo} 
\frac{\brack{ Mf,w } + \brack{ Au,w }}{\norm {(f,u) }_{\aLtwo \times \aHtwo}  \norm{ w }_{\aiLtwo}} \geq 1,
\end{equation*}
see \eqref{B6} and \eqref{B9}. 

Recall that we decided to write our state equation \eqref{A2}-\eqref{A3}
on the non-standard variational form \eqref{B1}.  Throughout this paper
we assume that problem \eqref{A2}-\eqref{A3} admits a unique solution $u \in \bHtwo$ for
every $f \in \Ltwo$, and that
\begin{equation}
\label{B2}
\norm{ u }_{\Htwo} \leq c_1 \norm{ f }_{\Ltwo}. 
\end{equation}
This assumption is valid if $\Omega$ is convex or if $\Omega$ has a
$C^2$ boundary, see e.g. \cite{grisvard1985, BHac92}. Inequality
\eqref{B2} is a key ingredient of the proof of our next lemma.
%% This results is presented in ``Satz 9.1.15'' (page 196) and ``Satz 9.1.17''
%% on pages 196 and 199, respectively, in the German edition of the book. 

\begin{lemma}  
\label{pre_coercivity}
There exists a constant $c_2$, which is independent of $\alpha > 0$, such that 
\begin{equation*}
  \begin{aligned}
    \alpha \norm{ f }_{\Ltwo}^2 + \norm{ u }_{\LtwoB}^2
    & \geq c_2 \left(
      \alpha \norm{ f }_{\Ltwo}^2 
      + \alpha \norm{ u }_{\Htwo}^2 
      + \norm{ u }_{\LtwoB}^2 
    \right) 
    \\
    & =  c_2 \norm{ (f,u) }_{\aLtwo \times \aHtwo}^2 
  \end{aligned}
\end{equation*}
for all $(f,u) \in \Ltwo\times \bHtwo$ such that
\begin{equation}
\label{B10}
(f,\phi)_{\Ltwo} + (-\Delta u + u,\phi)_{\Ltwo}  = 0 \quad \forall \phi \in \Ltwo.
\end{equation}
\end{lemma}
\subsubsection*{Proof} 
If $(f,u)$ satisfies \eqref{B10}, then 
\begin{equation*}
  \norm{ u }_{\Htwo} \leq c_1 \norm{ f }_{\Ltwo}, 
\end{equation*}
see the discussion of \eqref{B2}. Let $\theta = (1+ c_1^2)^{-1} \in (0,1)$, and it follows that
\begin{equation*}
  \begin{aligned}
    \alpha \norm{ f }_{\Ltwo}^2 + \norm{ u }_{\LtwoB}^2
    &\geq 
    \alpha \theta  \norm{ f }_{\Ltwo}^2 
    +\alpha\frac{1-\theta}{c_1^2} \norm{ u }_{\Htwo}^2 
    + \norm{ u }_{\LtwoB}^2 \\ 
    &\geq
      \frac{1}{1+ c_1^2}
    \left(
      \alpha \norm{ f }_{\Ltwo}^2 
      + \alpha \norm{ u }_{\Htwo}^2 + \norm{ u }_{\LtwoB}^2 
    \right).
  \end{aligned}
\end{equation*}
\rule{2mm}{2mm} \\
This result may also be written in the form 
\begin{equation*}
  \begin{split}
    \brack*{ 
      \left[
        \begin{matrix}
          \alpha M & \ze \\
          \ze & M_{\partial}  
        \end{matrix}
      \right]
      \left[
        \begin{matrix}
          f \\ u 
        \end{matrix}
      \right]
      , 
      \left[ \begin{matrix}
          f \\ u 
        \end{matrix}\right]
    }
    &\geq
    c_2 \norm{ (f,u) }_{\aLtwo \times \aHtwo}^2 \\ 
  \end{split}
\end{equation*}
for all $(f,u) \in \aLtwo \times \aHtwo$ satisfying 
\begin{equation*}
  Mf+Au=0, 
\end{equation*}
where $M$, $M_{\partial}$ and $A$ are the operators defined in \eqref{B6}, \eqref{B7} and \eqref{B9}, respectively. 
\subsection{Boundedness}

Having established that the Brezzi conditions hold, with constants that are 
independent of $\alpha$, we next explore the boundedness of $\Aca$. 

\begin{lemma}
\label{bounded_1}
\begin{equation*}
\abs*{
\brack*{
\left[
\begin{matrix}
\alpha M & \ze \\
\ze & M_{\partial}  
\end{matrix}
\right]
\left[
\begin{matrix}
f \\ u 
\end{matrix}
\right]
, 
\left[
\begin{matrix}
\psi \\ \phi 
\end{matrix}
\right]
}
}
\leq
\sqrt{2} \norm{ (f,u) }_{\aLtwo \times \aHtwo} \norm{ (\psi,\phi) }_{\aLtwo \times \aHtwo}   
\end{equation*}
for all $(f,u),(\psi, \phi) \in \aLtwo \times \aHtwo$. 
\end{lemma}
\subsubsection*{Proof} 
Recall the definitions \eqref{B6} and \eqref{B7} of $M$ and $M_{\partial}$, respectively. 
Since 
\begin{equation*}
\norm{ (f,u) }_{\aLtwo \times \aHtwo} = 
\sqrt{\alpha \norm{ f }_{\Ltwo}^2 + \alpha \norm{ u }_{\Htwo}^2 + \norm{ u }_{\LtwoB}^2}, 
\end{equation*}
we find, by employing the Cauchy-Schwarz inequality, that  
\begin{equation*}
    \begin{split}
\abs*{
\brack*{
\left[
\begin{matrix}
    \alpha M & \ze \\
    \ze & M_{\partial}  
  \end{matrix}
\right]
\left[
  \begin{matrix}
    f \\ u 
  \end{matrix}
\right]
, 
\left[
  \begin{matrix}
    \psi \\ \phi 
  \end{matrix}
\right]
}
} 
&=
\left| \alpha (f,\psi)_{\Ltwo} + (u,\phi)_{\LtwoB} \right| \\ 
&\leq \norm{ f }_{\aLtwo}  \norm{ \psi }_{\aLtwo} + \norm{ u }_{\LtwoB} \norm{ \phi }_{\LtwoB} \\ 
&\leq \sqrt{2}\sqrt{\norm{ f }_{\aLtwo}^2  \norm{ \psi }_{\aLtwo}^2 
+  \norm{ u }_{\LtwoB}^2  \norm{ \phi }_{\LtwoB}^2} \\
&\leq \sqrt{2} \norm{ (f,u) }_{\aLtwo \times \aHtwo} \norm{ (\psi,\phi) }_{\aLtwo \times \aHtwo}.  
  \end{split}
\end{equation*}
\rule{2mm}{2mm}

\begin{lemma}
\label{bounded_2}
\begin{equation*}
\abs*{
\brack*{
\left[
M \, \, A
\right]
\left[
\begin{matrix}
f \\ u 
\end{matrix}
\right]
, 
w
} 
}
\leq
\sqrt{3} \norm{(f,u) }_{\aLtwo \times \aHtwo} \norm{ w }_{\aiLtwo}   
\end{equation*}
for all $(f,u) \in \aLtwo \times \aHtwo$, $w \in \aiLtwo$. 
\end{lemma}
\subsubsection*{Proof} 
Again, we note that 
\begin{align*}
\norm{ (f,u) }_{\aLtwo \times \aHtwo} &= 
\sqrt{\alpha \norm{ f }_{\Ltwo}^2 + \alpha \norm{ u }_{\Htwo}^2 + \norm{ u }_{\LtwoB}^2}, \\
\norm{ w }_{\aiLtwo} &= \frac{1}{\sqrt{\alpha}} \norm{ w }_{\Ltwo}. 
\end{align*}
From the definitions of $M$ and $A$, see \eqref{B6} and \eqref{B9}, and
the Cauchy-Schwarz inequality, it follows that
\begin{equation*}
  \begin{split}
    \abs*{
    \brack*{
      \left[
        M \, \, A
      \right]
      \left[
        \begin{matrix}
          f \\ u 
        \end{matrix}
      \right]
      , 
      w
    }
    }
&= 
\abs*{
\brack*{
Mf,w
}
+ 
\brack*{
Au,w
}
} \\
&= \abs*{
(f, w)_{\Ltwo}
+ (-\Delta u + u, w)_{\Ltwo}  
} \\
&\leq  \left( \norm{ f }_{\aLtwo}  
+  \norm{ \laplace u }_{\aLtwo} 
+  \norm{ u }_{\aLtwo} \right)\norm{ w }_{\aiLtwo} \\
&\leq
\sqrt{3} \norm{ (f,u) }_{\aLtwo \times \aHtwo} \norm{ w }_{\aiLtwo}. 
  \end{split}
\end{equation*}
For the last equality, recall from \eqref{eq:H2_inner} that
$\norm{\laplace u}_\Ltwo = \norm{ \hess u}_\Ltwo \leq \norm{u}_\Htwo$
for all $u\in \bHtwo$.
\\ \rule{2mm}{2mm}

\subsection{Isomorphism}
We have verified that the Brezzi conditions hold, and that $\Aca$ is 
a bounded operator. Moreover, all constants appearing in the inequalities expressing 
these properties are independent of the regularization parameter $\alpha > 0$. 
Let 
\begin{align} 
\Vc  &= \aLtwo \times \aHtwo \times \aiLtwo, \\ 
\Vc ' &= \aLtwo' \times \aHtwo' \times \aiLtwo'.
\end{align} 

\begin{theorem}
\label{main}
The operator $\Aca$, defined in \eqref{B5}, 
%is an isomorphism taking $\Vc $ onto $\Vc '$, bounded uniformly for $\alpha>0$. That is, it 
is bounded and continuously invertible for $\alpha > 0$ in the sense that for all nonzero $x\in \Vc $,
\begin{equation} \label{Aa_continuous_estimates}
c \leq \sup_{ 0\neq y\in\Vc } \frac{\brack{\Aca x, y}}{\norm{y}_{\Vc }\norm{x}_{\Vc }} \leq  C,% \norm{x}_{\Vc }.
%\underline{c} \leq \norm{\Aca^{-1}}_{\Lc(\Vc ', \Vc )}^{-1}
%= \inf_{x\ne 0}\sup_{y\neq 0} \frac{\brack{A x, y}}{\norm{x} \norm{y}} \leq \sup_{x\ne 0}\sup_{y\neq 0} \frac{\brack{A x, y}}{\norm{x} \norm{y}} 
%< \norm{\Aca}_{\Lc(\Vc , \Vc ')} \leq \overline{c}.
% \| \Aca \|_{\Lc(\Vc , \Vc ')} \leq c_3 \quad\mbox{and} \quad \| \Aca^{-1} \|_{\Lc(\Vc ',\Vc )} \leq c_4, 
\end{equation}
for some positive constants $c$  and $C$  that are independent of $\alpha > 0$. In particular, 
\begin{equation*}
  \norm{\Aca^{-1}}_{\Lc(\Vc ',\Vc )} \leq c^{-1}  \quad\mbox{and} \quad \norm{\Aca}_{\Lc(\Vc , \Vc ')} \leq C.
\end{equation*}
\end{theorem}
\subsubsection*{Proof}
This result follows from Lemma \ref{pre_inf-sup}, Lemma \ref{pre_coercivity},
Lemma \ref{bounded_1}, Lemma \ref{bounded_2} and Brezzi theory for
saddle point problems, see \cite{brezzi}. \\
\rule{2mm}{2mm}

\subsection{Estimates for the discretized problem}
The stability properties \eqref{Aa_continuous_estimates} is not
necessarily inherited by discretizations.  However, the structure used
to prove the so-called ``inf-sup condition'' in Lemma
\ref{pre_inf-sup} is preserved in the discrete system provided that
 the same discretization is employed for the control and the Lagrange
multiplier.  Furthermore, the boundedness properties, Lemma
\ref{bounded_1} and Lemma \ref{bounded_2}, certainly also hold for
conforming discretizations. 

It remains to adress the coercivity condition, Lemma
\ref{pre_coercivity}, for the discretized problem.  We consider finite
dimensional subspaces $U_h\subset U = \bHtwo$ and $W_h\subset W =
\Ltwo$.  For certain choices of $U_h$ and $W_h$, the estimate of Lemma
\ref{pre_coercivity} carries over to the finite-dimensional setting.
\begin{lemma} \label{compatibility_lemma} 
  Assume $U_h\subset U$ and
  $W_h\subset W$, such that $(1-\laplace) U_h \subset W_h$.
  Then 
  \begin{equation}\label{discrete_coercivity}
    \alpha \norm{f_h}_\Ltwo^2 + \norm{u_h}_\LtwoB^2 \geq c_2 \norm{(f_h,u_h)}_{\aLtwo \times \aHtwo}^2
  \end{equation} 
  for all $(f_h,u_h) \in W_h\times U_h$ such that
  \begin{equation}\label{B10_discrete}
    (f_h, \phi_h)_\Ltwo + (u_h -\laplace u_h, \phi_h)_\Ltwo = 0 \quad \forall \phi_h \in W_h.
  \end{equation}
\end{lemma}
\subsection*{Proof}
Assume that $(1-\laplace) U_h \subset W_h$, and that
\eqref{B10_discrete} holds for $(f_h,u_h)\in W_h\times U_h$. Then $f_h
+(1-\laplace) u_h \in W_h$, and \eqref{B10_discrete} implies $f_h +(1-\laplace) u_h =
0$. Therefore, $(f_h, u_h)$ satisfies \eqref{B10} and the estimate
\eqref{discrete_coercivity} follows from Lemma
\ref{pre_coercivity}.  \\ \rule{2mm}{2mm}

If the discretization is chosen such that Lemma
\ref{compatibility_lemma} is satisfied, then the estimates
\eqref{Aa_continuous_estimates} carries over to discretized system.
More precisely, we have
\begin{align}
%\label{C11}
\norm{ \Acah }_{\Lc(\Vc _h,\Vc _h')} \leq \norm{ \Aca }_{\Lc(\Vc ,\Vc ')}, \quad\mbox{and}\quad
\norm{ \Acah^{-1} }_{\Lc(\Vc _h',\Vc _h)} \leq \norm{ \Aca^{-1} }_{\Lc(\Vc ',\Vc )},
\label{Aa_discrete_estimates}
\end{align}
where $\Vc _h = W_h\times U_h \times W_h \subset \Vc $, equipped with the
inner prdocut of $\Vc $, and $\Acah$ is discrete counterpart to $\Aca$,
defined by setting $\brack{\Acah x_h,y_h} = \brack{\Aca x_h, y_h}$ for
all $x_h,y_h\in \Vc _h$.

If the state is discretized with $C^1$-conforming bicubic
Bogner-Fox-Schmit rectangles, as in Section
\ref{Numerical_experiments}, then Lemma \ref{compatibility_lemma} is
satisfied if the control and Lagrange multiplier is discretized with
discontinuous bicubic elements on the same mesh.  For triangular
meshes, one could choose Argyris triangles for the state variable and
piecewise quintic polynomials for the control and Lagrange multiplier
variables.

We remark that Lemma \ref{compatibility_lemma} provides a sufficient,
but not necessary criterion for stability of the discrete problem, and
usually may imply far more degrees of freedom in the discrete space
$W_h\subset W$ than is actually needed.  The usefulness of Lemma
\ref{compatibility_lemma} is that the estimates
\eqref{Aa_discrete_estimates} can, in principle, always be obtained
by choosing a sufficiently large space for the control and Lagrange
multiplier.

\section{Preconditioning}\label{Preconditioning}
The linear problem \eqref{B5} is of the form
\begin{equation}
  \label{P:Axeqb}
  \Ac x = b.
\end{equation}
where $x$ is sought in a Hilbert space $\Vc $, the right hand side $b$ is
in the dual space $\Vc '$, and $\Ac $ is a self-adjoint continuous mapping
of $\Vc $ onto $\Vc '$. Iterative methods for linear problems are most
often formulated for operators mapping $\Vc $ into itself, and can not
be directly applied to the linear system \eqref{P:Axeqb}, as described
in \cite{Mar11}. If we want to apply such methods to \eqref{P:Axeqb},
then we need to introduce a continuous operator mapping $\Vc '$
isomorphically back onto $\Vc $. More precisely, if we have a continuous
operator
\begin{equation*}
  \Bc: \Vc ' \rightarrow \Vc , 
\end{equation*}
then $\Mc  = \Bc \Ac :\Vc \rightarrow\Vc $ is continuous and has the desired
mapping properties, and if $\Bc $ is an isomorphism, the solutions to
\eqref{P:Axeqb} coincides with the solutions to the problem
\begin{equation}
   \Mc  x = \Bc  \Ac  x =  \Bc  b. \label{P:BAxeqBb}
\end{equation}

In this paper we shall consider $\Bc  \in \Lc(\Vc ,\Vc ')$ a preconditioner 
if $\Bc $ is self-adjoint and positive definite. This implies that
$\Bc ^{-1}$ is self-adjoint and positive definite as well, and hence
$\Bc ^{-1}$ defines an inner product on $\Vc $ by setting
\begin{equation} \label{P:Riesz}
  \pp{x,y} = \brack{ \Bc ^{-1} x, y }, \quad\quad  x,y \in \Vc .
\end{equation}
This inner product has the crucial property of making $\Mc $ self-adjoint, in the sense that
\begin{equation}
  \pp{\Mc  x,y} = \brack{\Ac  x,y} =  \brack{\Ac  y,x} = \pp{\Mc  y,x}. \label{P:symmetric}
\end{equation}

Conversely, given any inner product on $\pp{\subs,\subs}$ on $\Vc $,
the Riesz-Fr\'echet theorem provides a self-adjoint positive definite
isomorphism $\Bc :\Vc '\rightarrow \Vc $ such that \eqref{P:Riesz} and
\eqref{P:symmetric} hold, and we say that $\Bc $ is the Riesz operator
induced by $\pp{\subs,\subs}$.  This establishes a one-to-one
correspondence between preconditioners and Riesz operators on $\Vc '$.
Since the Riesz operator is an isometric isomorphism, the operator
norm of $\Bc \Ac $ coincides with the operator norm of $\Ac $.  We formulate
this well-known fact here in a lemma for the sake of self-containedness. We refer to
\cite{Mar11, gunnel2014note} for a more in-depth discussion of
preconditioning and its relation to Riesz operators.
\begin{lemma} \label{precond:equivalence}
  Let $\Vc $ be a Hilbert space, and let $\Ac :\Vc \rightarrow \Vc '$ be a
  self-adjoint isomorphism, and assume that $\Bc $ is the Riesz operator
  induced by the inner product on $\Vc $, or
  equivalently, that the inner product on $\Vc $  is defined by the self-adjoint
  positive definite isomorphism $\Bc ^{-1}:\Vc \rightarrow \Vc '$.  Then $\Bc  \Ac  : \Vc  \rightarrow \Vc $ 
  is an isomorphism, self-adjoint in the inner product on $\Vc $,  with
  \begin{equation*}
    \norm{\Bc \Ac }_{\Lc(\Vc , \Vc )} = \norm{\Ac }_{\Lc(\Vc , \Vc ')} 
    \quad \mbox{and}\quad
    \norm{(\Bc \Ac )^{-1}}_{\Lc(\Vc , \Vc )} = \norm{\Ac ^{-1}}_{\Lc(\Vc ', \Vc )}. 
  \end{equation*}
  In particular, the condition number of $\Bc  \Ac $ is given by
  \begin{equation*}
    \cond(\Bc  \Ac ) = \norm{\Ac ^{-1}}_{\Lc(\Vc ', \Vc )}\norm{\Ac }_{\Lc(\Vc , \Vc ')}.
  \end{equation*}
\end{lemma}
\subsubsection*{Proof}
Since $\Ac $ is self-adjoint, $\Mc  = \Bc \Ac $ is
self-adjoint with respect to the inner product on $\Vc $.  From the
Riesz-Fr\'echet theorem we have $\norm{\Ac  x}_{\Vc '} = \norm{\Bc \Ac  x} =
\norm{\Mc  x}$, and we obtain following identity for the operator norm
of $\Mc $.
\begin{equation*}%   \label{C:M}
  \begin{split}
  \norm{\Mc }_{\Lc(\Vc , \Vc )} = \sup_{x\neq 0} \frac{ \norm{\Mc  x}_\Vc }{\norm{x}_\Vc } 
  &= \sup_{x\neq 0} \frac{ \norm{\Ac  x}_{\Vc '}}{\norm{x}_\Vc } \\
  &= \sup_{x\neq 0}\sup_{y\neq 0} \frac{\brack{\Ac  x, y}}{\norm{x}_\Vc \norm{y}_\Vc }
  = \norm{\Ac}_{\Lc(\Vc , \Vc ')}.
  \end{split}
\end{equation*}
A similar identity is obtained for the norm of the inverse operator,  
\begin{equation*}%\label{C:Minv}
  \begin{split}
      \norm{\Mc ^{-1}}_{\Lc(\Vc , \Vc )}  
  &= \sup_{x\neq 0} \frac{ \norm{\Mc ^{-1} x}_\Vc }{\norm{x}_\Vc } 
   \\ &= \left(\inf_{x\neq 0} \frac{ \norm{\Mc  x}_\Vc }{\norm{x}_\Vc } \right)^{-1}
%  = \inf_{x\neq 0}\sup_{y\neq 0} \frac{ \pp{\Mc  x, y}}{\norm{x}\norm{y}} % ^{-1}
  \\ &= \left(\inf_{x\neq 0}\sup_{y\neq 0} \frac{\brack{\Ac x, y}}{\norm{x}_\Vc \norm{y}_\Vc }\right)^{-1}
  = \norm{\Ac^{-1}}_{\Lc(\Vc ', \Vc )}.
  \end{split}
\end{equation*}
\rule{2mm}{2mm}
\\
We say that a preconditioner $\Bca$ for $\Aca$ is robust
with respect to the parameter $\alpha$ if $\cond(\Bca\Aca)$
is  bounded uniformly in $\alpha$. The significance of Lemma
\ref{precond:equivalence} is that such a robust preconditioner can be
found by identifying (parameter-dependent) norms in which $\Aca$ and
$\Aca^{-1}$ are both uniformly bounded.

\subsection{Parameter-robust minimum residual method}
In Section \ref{Analysis} stability of $\Aca$ was shown in the
$\alpha$-dependent norms defined in \eqref{B0.1}-\eqref{B0.3}.  The
preconditioner provided by Lemma \ref{precond:equivalence} is the
Riesz operator induced by the weighted norms. This operator $\Bca: \Vc '
\rightarrow \Vc $ takes the form
\begin{equation} \label{def_Bca}
\Bca=
\left[
\begin{matrix}
\alpha M & \ze & \ze \\
\ze & \alpha R +M_{\partial} & \ze \\
\ze &  \ze & \frac{1}{\alpha} M 
\end{matrix}
\right]^{-1}
\end{equation}
where $R:\bHtwo\rightarrow \bHtwo'$ is the operator induced by the
$\Htwo$ inner product, i.e.  $\brack{Ru,v} = (u,v)_{\Htwo}$.

Since $\Aca$ is self-adjoint, the
preconditioned operator $\Bca\Aca:\Vc \rightarrow\Vc $ is self-adjoint in the
inner product on $\Vc $.  Consequently we can apply the
minimum residual method (\MINRES) to the problem
\begin{equation*}
  \Bca \Aca x = \Bca b. 
\end{equation*}
\begin{theorem}
\label{number_of_iterations}
Let $\Aca$ be the operator defined in \eqref{B5} and $\Bca$ the operator
defined in \eqref{def_Bca}.  Then there exists an upper bound,
independent of $\alpha$, for the convergence rate of \MINRES{} applied to the preconditioned system
\begin{equation*}
  \Bca \Aca x = \Bca b.
\end{equation*}
In particular there exists an upper bound, independent of $\alpha$, for the number of iterations 
needed to reach the stopping criterion  \eqref{N3}.
\end{theorem} 
\subsubsection*{Proof}
A crude upper bound for the convergence rate (more precisely, the
two-step convergence rate) of \MINRES{} is given by
\begin{equation*}
  \norm{\Bca\Aca(x-x_{2m})}_{\Vc } \leq \left(\frac{1-\cond}{1+\cond}\right)^{m} \norm{\Bca\Aca(x-x_{0})}_{\Vc }
\end{equation*}
where $\cond = \cond(\Bca\Aca)$ is the condition number of $\Bca\Aca$,
see e.g. \cite{Mar11}.  From Lemma \ref{precond:equivalence} and
\eqref{Aa_continuous_estimates} we determine that $\cond$ is bounded
independently of $\alpha$, with
\begin{equation}
  \label{P:cond}
  \begin{split}
      \cond &=  \norm{(\Bca\Aca)^{-1}}_{\Lc(\Vc , \Vc )}  \norm{\Bca\Aca}_{\Lc(\Vc , \Vc )} 
       \\ &=  \norm{\Aca^{-1}}_{\Lc(\Vc ', \Vc )}  \norm{\Aca}_{\Lc(\Vc , \Vc ')}
       \\ &\leq c^{-1}C.
  \end{split}
\end{equation}
\rule{2mm}{2mm} \\ \\
In practical applications, the operator $\Bca$ will be replaced with a
less computationally expensive approximation $\widehat\Bca$.  Ideally
$\widehat\Bca$ will be spectrally equivalent to $\Bca$, in the sense
that the condition number of $\widehat\Bca \Bca^{-1}$ is bounded,
independently of $\alpha$. Then the preconditioned system reads
\begin{equation*}
  \widehat\Bca \Aca x = \widehat\Bca b,
\end{equation*}
and the upper bound for the convergence rate is determined by the conditioned number 
$\cond(\widehat\Bca\Aca) \leq \cond(\widehat\Bca\Bca^{-1})\cond(\Bca\Aca^{-1})$.

\subsubsection*{Remark}
In this paper we only consider the minimum residual method, and we
therefore require that the preconditioner is self-adjoint and positive
definite.  More generally, if other Krylov subspace methods are
to be applied to \eqref{B5}, then preconditioners lacking symmetry or
definiteness may be considered.

We mention in particular that a preconditioned conjugate gradient
method for problems similar to \eqref{B5} was proposed in \cite{s-z},
based on a clever choice of inner product.

% not necessary to deman SPD preconditioner
% other Krylov methods 

\section{Generalization}
\label{Generalization} 
Is our technique applicable to other problems than \eqref{A1}-\eqref{A3}? 
We will now briefly explore this issue, and show that the preconditioning scheme  
derived above yields $\alpha$ robust methods for a class of problems. 

The scaling \eqref{B0.1}-\eqref{B0.3} was also investigated in
\cite{Nie10}, but for a family of abstract problems posed in terms of
Hilbert spaces. More specifically, for general PDE-constrained
optimization problems, subject to Tikhonov regularization, and with
linear state equations. But in \cite{Nie10} no assumptions about the
control, state or observation spaces were made, except that they were
Hilbert spaces. Under these circumstances, it was proved that the
coercivity and the boundedness, of the operator associated with the
KKT system, hold with $\alpha$-independent constants. Nevertheless, in
this general setting, the inf-sup condition involved an
$\alpha$-dependent constant, which, eventually, yielded theoretical
iteration bounds of order $O([\log\left( \alpha^{-1} \right)]^2)$ 
for \MINRES{}.

In the present paper we were able to prove an $\alpha$-robust 
inf-sup condition for the model problem \eqref{A1}-\eqref{A3}.  This is possible because both the
control $f$ and the dual/Lagrange-multiplier $w$ belong to $\Ltwo$.
%This is possible assuming that both the
%control $f$ and the dual/Lagrange-multiplier $w$ belong to $\Ltwo$, but leads to  
%extra regularity requirement for the state.  
From a more general perspective, it turns out that this is the
property that must be fulfilled in order for our approach to be
successful: The control space and the dual space, associated with the
state equation, must coincide.  This will usually lead to additional
regularity requirements for the state space.

Motivated by this discussion, let us consider an abstract problem of the form: 
\begin{equation}
\label{G1}
\min_{f \in W, \, u \in U} 
\left\{ 
\frac{1}{2} \norm{ Tu - d  }^2_O 
+ \frac{1}{2} \alpha \norm{ f }^2_W
\right\}
\end{equation}
subject to 
\begin{equation}
\label{G2}
\brack{ Au, w } +  (f,w)_W =0, 
\quad \forall w \in W. 
\end{equation}
Here, 
\begin{itemize} 
\item $W$ is the dual \underline{and} control space, 
\item $U$ is the state space, 
\item $O$ is the observation space, 
\item $W$, $U$ and $O$ are Hilbert spaces. 
\end{itemize}
Let us assume that 
\begin{itemize}
\labeleditem{(A1)}{G:A1}
  $A:U \rightarrow W'$ is a continuous linear operator with
  closed range.  In particular, there is a constant $c_1$ such that for all $u \in U$, 
  \begin{equation*}
    \norm{u}_{U/ \ker A} = \inf_{\tilde u \in \ker A} \norm{u-\tilde u}_U \leq c_1\norm{A u}_{W'}.
  \end{equation*}
%
%\item $A^{-1}:W' \rightarrow U$ is bounded, 
%
\labeleditem{(A2)}{G:A2}
  $T:U \rightarrow O$ is linear and bounded, and invertible on the kernel of $A$.
  That is, there is a constant $c_2$ such that for all $u\in \ker A$,
  \begin{equation*}
    \norm{u}_U \leq c_2 \norm{T u}_O. 
  \end{equation*}
\end{itemize}
It then follows that the KKT system associated with 
\eqref{G1}-\eqref{G2} is well-posed for every $\alpha > 0$: 
Determine $(f,u,w) \in W \times U \times W$ such that  
\begin{equation}
\label{G3}
\underbrace{\left[
\begin{matrix}
\alpha M & \ze & M' \\
\ze & K & A' \\
M &  A & 0 
\end{matrix}
\right]}_{=\Aca}
\left[
\begin{matrix}
f \\ u \\ w
\end{matrix}
\right]
=
\left[
\begin{matrix}
0 \\ 
\tilde K d \\ 
0
\end{matrix}
\right], 
\end{equation}
where 
\begin{align}
\label{G4}
 M&: W \rightarrow W', &f &\mapsto (f,\subs)_W, \\ 
\label{G5}
 K&: U \rightarrow U', &u &\mapsto (Tu,T\subs)_O, \\ 
 \tilde K&: O \rightarrow U', &d &\mapsto (d,T\subs)_O,
%
%\label{G6}
% Q&: O \rightarrow U', &d &\mapsto (d,T \subs)_O.
\end{align} 
Note that, compared with \eqref{N1}, the boundary observation matrix $M_{\partial}$ has been 
replaced with the general observation operator $K$ in \eqref{G3}.  

We introduce scaled norms as follows.
\begin{align*}
 \norm{ f }_{W_{\alpha}}^2 &= \alpha \norm{ f }_W^2, \\
 \norm{ u }_{U_{\alpha}}^2 &= \alpha \norm{Au}_{W'}^2 + \norm{ Tu }_O^2, \\
 \norm{ w }_{W_{\alpha^{-1}}}^2 &= \frac{1}{\alpha} \norm{ w }_W^2.
\end{align*}
We first show that $\norm{\subs}_{U_\alpha}$ is indeed a norm on $U$ when assumptions \ref{G:A1} and \ref{G:A2} hold.
It suffices to show that $\norm{\subs}_{U_\alpha}$ is a norm equivalent to
$\norm{\subs}_U$ when $\alpha=1$. We have
\begin{equation}\label{G:norm_euiv1}
  \norm{Tu}_O + \norm{Au}_{W'} 
  \leq \big(\norm{T}_{\Lc( U,O)}  + \norm{A}_{\Lc(U,W')}\big) \norm{u}_U,
\end{equation}
and letting $\pi$ denote the orthogonal projection of $U$ onto $\ker A$,
\begin{equation}
  \begin{split}
    \norm{u}_U&\leq \norm{\pi u}_U + \norm{u-\pi u}_U
    \\ &\leq c_2 \norm{T\pi u}_O + \norm{u-\pi u}_U
    \\ &\leq c_2 \norm{T u}_O + \big(1+c_2 \norm{T}_{\Lc(U,O)}\big)\norm{u-\pi u}_U
%    = c_1 \norm{T u}_O + (1+c_1 c_2)\norm{u}_{U/ \ker A}
    \\ & \leq c_2 \norm{T u}_O +c_1\big(1+c_2\norm{T}_{\Lc(U,O)}\big) \norm{Au}_{W'}.
  \end{split}
\end{equation}
Here the last inequality follows from $\norm{u-\pi u}_{U} =
\inf_{\tilde u\in \ker A}\norm{u-\tilde u}_U$ and assumption \ref{G:A1}.

We set $\Vc  = W_\alpha \times U_{\alpha} \times W_{\alpha^{-1}}$. As in
Section \ref{Analysis}, $\Aca:\Vc  \rightarrow \Vc '$ can be shown to be
an isomorphism, with parameter-independent estimates obtained in the
weighted norms.
\begin{theorem}
\label{GG:isomorphism}
  There exists positive constants $c$ and $C$,
  independent of $\alpha$, such that for all
  nonzero $x \in \Vc $,
  \begin{equation} \label{G:isomorphism}
    c \leq \sup_{0\neq y \in \Vc } \frac{\brack{\Aca x, y}}{ \norm{x}_{\Vc }  \norm{y}_{\Vc }} \leq  C.
  \end{equation}
\end{theorem}
We omit the full proof, which is analogous to that of Theorem
\ref{main}. The crucial part is the ``inf-sup condition'' of Lemma \ref{pre_inf-sup}, which is
easily shown to hold in the abstract setting:
\begin{equation*} 
  \begin{aligned}
    \sup_{(f,u) \in W_\alpha \times U_\alpha} 
    \frac{(f,w)_{W} + \brack{A u,w}}{\norm{ (f,u) }_{W_\alpha\times U_\alpha}  \norm{ w }_{W_{\alpha^{-1}}}} &
    \geq
    \frac{(w,w)_{W}}{\norm{ (w,0) }_{W_\alpha \times U_\alpha} \norm{ w }_{W_{\alpha^{-1}}}}
    %\frac{\| w \|^2_{\Ltwo}}{\sqrt{\alpha}\| w \|_{\Ltwo} \, (\sqrt{\alpha})^{-1} \| w \|_{\Ltwo}} 
    = 1.  
  \end{aligned}
\end{equation*}
The coercivity condition of Lemma \ref{pre_coercivity} naturally holds
in the prescribed norm on $U_\alpha$, since for $(f,u)\in W\times U$
such that $Au = Mf$,
\begin{equation*}
  \alpha \norm{f}_W^2 + \norm{T u}_O^2 
  = \frac{\alpha}{2} \norm{f}_W^2 + \frac{\alpha}{2} \norm{Au}_{W'}^2 +  \norm{T u}_O^2  
  \geq \frac{1}{2} \left(\norm{f}_{W_\alpha}^2 + \norm{u}_{U_\alpha}^2 \right).
\end{equation*}

Note that the weighted norm now depends on $A$, and as consequence,
the estimates become $A$-independent. In fact, we obtain bounds for
the constants $c$ and $C$ which are independent of $\alpha$ as well as the  
operators appearing in \eqref{G1}-\eqref{G2}. This is postponed to
the next section, where sharp estimates are obtained for
\eqref{G:isomorphism}.

With the estimates \eqref{G:isomorphism}, Lemma
\ref{precond:equivalence} provides a preconditioner for the operator
$\Aca$, given as
\begin{equation}\label{G:precond}
\Bca=
\left[
\begin{matrix}
\alpha M & \ze                   & \ze \\
\ze      & \alpha A' M^{-1} A + K & \ze \\
\ze      & \ze                   & \frac{1}{\alpha} M 
\end{matrix}
\right]^{-1} .   
\end{equation}
The condition number of $\Bca \Aca$ will be bounded independently of
$\alpha$.  It is, however, not clear how to find a computationally
efficient approximation of $\Bca$ in the abstract setting of
\eqref{G1}-\eqref{G2}.

\subsubsection*{Example 1}
The problem \eqref{A1}-\eqref{A3} fits in the abstract framework
presented in this section when we assume that the state has $\Htwo$
regularity. We set $W= \Ltwo$, $U=\bHtwo$, $A = 1-\Delta$, and
$T:\bHtwo\rightarrow \LtwoB$ is a trace operator, see (\ref{G5}). Since $A$ is a continuous
isomorphism, assumptions \ref{G:A1} and \ref{G:A2} are both valid.  The inner
product on $U_\alpha$ takes the form
\begin{equation*}
  \begin{split}
     (u, v)_{U_\alpha} &= \brack{K u, v} + \alpha \brack{A M^{-1} A u, v}
     \\ &= \int_{\partial\Omega} u v \, ds + \alpha\int_\Omega (u-\laplace u)(v-\laplace v) \, dx
     \\ &= \int_{\partial\Omega} u v \, ds + \alpha\int_\Omega D^2u : D^2v + 2 \grad u \cdot \grad v \, +uv \,dx,
  \end{split}
\end{equation*}
where $D^2 u$ denotes the Hessian of $u$, and the last equality
follows from the boundary condition $\partial u /\partial \mathbf{n} =
0$ imposed on $\bHtwo$. The resulting preconditioner is the one that
was used in the numerical experiments, detailed in Section
\ref{Numerical_experiments}, and it is spectrally equivalent to the
preconditioner defined in \eqref{def_Bca}.

\subsubsection*{Example 2}
Let $U$, $W$, and $K$ be as in Example 1, but let us
set $A = -\laplace$.  Now $A$ has non-trivial kernel, consisting of the
a.e. constant functions, and for constant $u$ we  have
\begin{equation*}
  \norm{Tu}_\LtwoB =  \sqrt{\frac{\abs{\partial \Omega}}{\abs{\Omega}}} \norm{u}_\bHtwo.
\end{equation*}
Since assumptions \ref{G:A1} and \ref{G:A2} are valid, the optimality system
is still well-posed.  In this case the inner product on $U_\alpha$ is
given by
\begin{equation*}
  \begin{split}
     (u, v)_{U_\alpha} = \int_{\partial\Omega} u v \, ds + \alpha\int_\Omega D^2u:D^2v \,dx.
  \end{split}
\end{equation*}

\subsubsection*{Example 3}
Let us consider the ``prototype'' problem: 
\begin{equation*}
\min_{f, \, u} 
\left\{ 
  \frac{1}{2}\norm{ u - d }_{\Ltwo}^2
  + \frac{\alpha}{2}\norm{ f }_{\Ltwo}^2
\right\}
\end{equation*}
subject to 
\begin{alignat*}{2}
- \laplace u + u + f&= 0  &&\quad \mbox{in } \Omega, \\
\frac{\partial u}{\partial \mathbf{n}} &= 0  &&\quad \mbox{on } \partial \Omega. 
\end{alignat*}
Note that we here consider the case in which observation data is assumed
to be available throughout the entire domain $\Omega$ of the state
equation. 

If the usual variational form of the PDE is used, i.e.,
\begin{equation}\label{G:ex_var1}
(u,w)_{\Hone} + (f,w)_{\Ltwo} = 0, \quad\forall w \in \Hone, 
\end{equation}
then the control space equals $\Ltwo$, whereas the dual
space is $\Hone$.  The preconditioning strategy presented in this
section is therefore not applicable.

If instead we can assume $\Htwo$-regularity, we can use the variational form
\begin{equation}\label{G:ex_var2}
    (- \laplace u +u,w)_{\Ltwo} + (f,w)_{\Ltwo} =0, \quad \forall w \in \Ltwo.
\end{equation}
Now, the control and dual spaces both equal $\Ltwo$.  The methodology
presented in this section can thus be applied, and a robust
preconditioner is obtained. Compared with the preconditioner for the
problem with boundary observations only, see Section
\ref{Preconditioning}, equation (\ref{def_Bca}), the only change is
the replacement of $M_{\partial}$, in the $(2,2)$ block of $\Bca$ with
$M$. 

We remark that in \cite{s-z} and \cite{Pea12}, parameter-robust
preconditioners were proposed for the ``prototype'' problem, using the
standard variational formulation \eqref{G:ex_var1} of the PDE. Those
methods do not require improved regularity for the state
space. Instead, they require that observations are available
throughout the computational domain.

\section{Eigenvalue analysis}
\label{Alternative}
In Section \ref{Generalization} it was shown that the condition
number of $\Bca\Aca$, with $\Aca$ defined in \eqref{G3} and $\Bca$
defined in \eqref{G:precond}, can be bounded independently of
$\alpha$, as well as independently of the operators appearing in
\eqref{G1}-\eqref{G2}.  Moreover, the numerical experiments indicate 
that the eigenvalues are contained in three intervals, independently of
the regularization parameter $\alpha$, see Figure \ref{fig:1}. In this
section we detail the structure of the spectrum of the preconditioned
system considered in Section \ref{Generalization}, and we obtain sharp
estimates for the constants appearing in Theorem \ref{GG:isomorphism}.

We consider self-adjoint linear operators
$\mathcal{A}_{\alpha}$ and $\mathcal{B}_\alpha $,
\begin{equation} \label{aa:systems}
  \mathcal{A}_{\alpha} = 
  \left[ \begin{matrix}    \alpha  M &   \ze          &  M'   \\
                              \ze    &   K            &  A '  \\
                               M     &   A            &   0   \end{matrix} \right]
\quad\mbox{and}\quad
  \mathcal{B}_\alpha^{-1}  = 
  \left[ \begin{matrix}  \alpha  M &    \ze        &   \ze \\
                            \ze    & K  + \alpha R &   \ze \\
                            \ze    &    \ze        &   \alpha^{-1} M      \end{matrix} \right]
\end{equation}
where $R$ is defined by
\begin{equation}
  R =  A '  M ^{-1} A  . \label{aa:def_R}
\end{equation}
We assume that $A:U\rightarrow W'$ and $M:W\rightarrow W'$ are
continuous operators, for some Hilbert spaces $U$ and $W$. In addition
we will make use of the following assumptions.
\begin{itemize}
\labeleditem{(B1)}{AA:A1}
  $ M $ is a self-adjoint and positive definite,
\labeleditem{(B2)}{AA:A2}
 $ K  + R$ is positive definite,
\labeleditem{(B3)}{AA:A3}
  $ K $ is self-adjoint and positive semi-definite.
\end{itemize}
Assumptions \ref{AA:A1}-\ref{AA:A3} ensure that $\mathcal{B}_\alpha $ is a
self-adjoint and positive definite.  In particular, assumptions \ref{AA:A1}-\ref{AA:A3} hold
for $\Aca$ as in \eqref{G3}, provided that the assumptions of Section
\ref{Generalization} hold.  For simplicity, we also assume that that
$\mathcal{A}_{\alpha}$ and $\mathcal{B}_\alpha $ are
finite-dimensional operators.
\begin{theorem} \label{aa:th}
  Let $p$, $q$, and $r$ be the polynomials
  \begin{equation*}
    p(\lambda) = 1-\lambda, \quad
    q(\lambda) = 1+\lambda p(\lambda), \quad
    r(\lambda) = p - \lambda q(\lambda).
  \end{equation*}
  Let $q_1<q_2$ and $r_1 < r_2 <r_3$ be the roots of $q$ and $r$,
  respectively.  The spectrum of $\mathcal{B}_\alpha  \mathcal{A}_{\alpha}$ is contained
  within three intervals, determined by the roots of $p$ and $r$, 
  independently of $\alpha$: 
  \begin{equation}
    \operatorname{sp}(\mathcal{B}_\alpha \mathcal{A}_{\alpha}) \subset [r_1,q_1] \cup [r_2,1] \cup [q_2,r_3].
    \label{aa:spectral_bounds}
  \end{equation}
  Consequently, the spectral condition number of $\mathcal{B}_\alpha  \mathcal{A}_{\alpha}$ is bounded, uniformly in $\alpha$,
  \begin{equation}
    k (\mathcal{B}_\alpha \mathcal{A}_{\alpha}) \leq \frac{r_3}{r_2} \approx 4.089.
    \label{aa:condition}
  \end{equation}
  If $ K $ has a nontrivial kernel, inequality \eqref{aa:condition}
  becomes an equality.
\end{theorem}

\subsubsection*{Proof}
  Consider the equivalent generalized eigenvalue problem
  \begin{align}
  \left[ \begin{array}{ccc}    
      \alpha  M &   \ze           &  M'   \\
        \ze     &   K             &  A '  \\
          M     &   A             &   0   \end{array} \right]
\left[\begin{matrix}  f  \\  u \\  w  \end{matrix} \right]
\quad = \quad \lambda\, 
  \left[ \begin{array}{ccc}    \alpha  M &    \ze     & \ze    \\
                           \ze        &   K  + \alpha R           &  \ze  \\
                           \ze    &   \ze          &   \alpha^{-1} M      \end{array} \right]
\left[\begin{matrix}  f  \\  u \\  w  \end{matrix} \right]
  \label{aa:gep}
\end{align}
We show that \eqref{aa:gep} admits no nontrivial solutions 
unless $\lambda$ is as in \eqref{aa:spectral_bounds}.
 
Since $ M $ is a self-adjoint isomorphism, by assumption \ref{AA:A1}, we can
rewrite \eqref{aa:gep} as the three identities
\begin{align}
  \alpha p  f   + w           & = 0,    \label{aa:gep1}\\
  p K   u     + A '  w   - \lambda \alpha R  u   & = 0,  \label{aa:gep2}\\
    f    +    M  ^{-1} A   u -\lambda \alpha^{-1}  w      & = 0.  \label{aa:gep3}         
\end{align}
Assume that $\lambda$ is not contained within the three closed intervals
of \eqref{aa:spectral_bounds}.
Then $p \neq 0$, and we can use  \eqref{aa:gep1} 
to eliminate $ f $ from \eqref{aa:gep3}.
\begin{equation}
\label{special1}
  \begin{split}
      0 &=  \alpha p ( f  +  M ^{-1}  A   u - \lambda \alpha^{-1} w ) 
        =  \alpha p  M ^{-1}  A   u  - (1 +\lambda p) w  \\
        &= \alpha p  M ^{-1}  A   u  -  q  w . 
  \end{split}
\end{equation}
Since $q$ is nonzero, we can use \eqref{special1} to eliminate $ w $ from \eqref{aa:gep2},
\begin{equation}
  \begin{split}
      0 &=  q (p  K   u  + A '  w   - \lambda \alpha R  u ) 
         = qp  K   u + \alpha(p - \lambda q)R  u \\
        &= qp  K   u + r R u,
        \label{aa:x2}
  \end{split}
\end{equation}
where the identity \eqref{aa:def_R} was used.  By assumption, $pq$ and
$r$ are both nonzero. Moreover, it can be easily seen that $pq$ and
$r$ have the same sign outside of the bounded intervals of
\eqref{aa:spectral_bounds}. From assumptions \ref{AA:A1}-\ref{AA:A3}, we conclude that
$qp K + r R$ is a self-adjoint definite operator.  Then \eqref{aa:x2}
only admits trivial solutions, hence  $\lambda$ can not be
an eigenvalue of $\mathcal{B}_\alpha \mathcal{A}_{\alpha}$.

The estimate \eqref{aa:condition} follows from
\eqref{aa:spectral_bounds}, noting that $\vert \operatorname{sp}(\mathcal{B}_\alpha
\mathcal{A}_{\alpha}) \vert \subset [r_2, r_3]$.  From \eqref{aa:x2} it
can be seen that the roots of $r$ are eigenvalues of
$\mathcal{B}_\alpha \mathcal{A}_{\alpha}$ if $\ker K $ is nontrivial.
\\
\rule{2mm}{2mm}

\subsection*{Remark}
If $A = (1 - \laplace):\bHtwo\rightarrow\Ltwo'$, then $R = A' M^{-1} A$
is characterized by a bilinear form $b(\cdot,\cdot)$ as in \eqref{N2.1}:
\begin{equation*}
  \begin{split}
   \brack[\big]{ A' M^{-1} A u, v } &=
   \int_\Omega   \laplace u \laplace v 
        +  2 \grad u \cdot \grad v 
        +    u v \dx
        \\
  & =  (u,v)_{\Htwo} + \int_\Omega   \grad u \cdot \grad v \dx = b(u, v)
  \end{split}
\end{equation*}
For discretizations $U_h\subset U$ and $W_h\subset W$ of $A$ such that
$A(U_h) \subset M(W_h)$, the discretization of $b$ coincides with
$A_h'M_h^{-1} A_h$. This follows from an argument similar to that in
the proof of Lemma \ref{compatibility_lemma}, and as a consequence,
Theorem \ref{aa:th} can be applied to the preconditioned discrete
systems considered in Section \ref{Numerical_experiments}.

% example
\section{Discussion}
\label{Conclusions}
Previously, parameter robust preconditioners for PDE-constrained
optimization problems have been successfully developed, provided that
observation data is available throughout the entire domain of the
state equation.  For many important inverse problems, arising in
industry and science, this is an unrealistic requirement. On the
contrary, observation data will typically only be available in
subregions, of the domain of the state variable, or at the boundary of
this domain. We have therefore explored the
possibility for also constructing robust preconditioners for
PDE-constrained optimization problems with limited observation data.

For an elliptic control problem, with boundary observations only, we
have developed a regularization robust preconditioner for the
associated KKT system.  Consequently, the number of \MINRES{} iterations
required to solve the problem is bounded independently of both
regularization parameter $\alpha$ and the mesh size $h$. In order to
achieve this, it was necessary to write the elliptic state equation on
a non-standard, and non-self-adjoint, variational form. If this
approach is employed, then the control and the Lagrange multiplier
will belong to the same Hilbert space, which leads
to extra regularity requirements for the state.  This fact makes it possible to
construct parameter weighted metrics such that the constants
appearing in the Brezzi conditions, as well as the
constants in the inequalities expressing the boundedness of the KKT
system, are independent of $\alpha$ and $h$. Consequently, the
spectrum of the preconditioned KKT system is uniformly bounded with
respect to $\alpha$ and $h$, which is ideal for the \MINRES{} scheme. 
These properties were illuminated through a series of numerical experiments, 
and the preconditioned \MINRES{} scheme handled our model problem excellently.   

The use of a non-self-adjoint form of the elliptic state equation
leads to additional challenges for constructing discretization schemes
and suitable multigrid methods.  More specifically, it becomes
necessary to implement a FE space approximating $H^2$. We accomplished
this by a $C^1$ discretization that is conforming in $H^2$. The method
employed does, however, have strong restrictions on the mesh, which
seemingly must be composed of rectangles.

\bibliographystyle{plain}
\bibliography{robust_precond_kkt_final}

\end{document}